\documentstyle[12pt]{article}
\begin{document}
\pagestyle{empty}
\def\lsim{\raise0.3ex\hbox{$<$\kern-0.75em\raise-1.1ex\hbox{$\sim$}}}
\def\gsim{\raise0.3ex\hbox{$>$\kern-0.75em\raise-1.1ex\hbox{$\sim$}}}
\def\noi{\noindent}
\def\nn{\nonumber}
\def\bea{\begin{eqnarray}}  \def\eea{\end{eqnarray}}
\def\beq{\begin{equation}}   \def\eeq{\end{equation}}
\def\sq{\hbox {\rlap{$\sqcap$}$\sqcup$}}
\def\beeq{\begin{eqnarray}} \def\eeeq{\end{eqnarray}}
\newcommand\mysection{\setcounter{equation}{0}\section}
\renewcommand{\theequation}{\thesection.\arabic{equation}}
\newcounter{hran} \renewcommand{\thehran}{\thesection.\arabic{hran}}
\renewcommand{\thefootnote}{\fnsymbol{footnote}}
\def\R{ {\rm R \kern -.31cm I \kern .15cm}}
\def\C{ {\rm C \kern -.15cm \vrule width.5pt \kern .12cm}}
\def\Z{ {\rm Z \kern -.27cm \angle \kern .02cm}}
\def\N{ {\rm N \kern -.26cm \vrule width.4pt \kern .10cm}}
\def\1{{\rm 1\mskip-4.5mu l} }

\vbox to 2 truecm {}
\centerline{\Large \bf Scattering Theory in the Energy Space}
\vskip 3 truemm  
\centerline{\Large \bf for a Class of Hartree Equations*}
\vskip 8 truemm
\centerline{\it Dedicated to Professor Walter A. Strauss on his 60th birthday} 

\vskip 1 truecm
\centerline{\bf J. Ginibre}
\centerline{Laboratoire de Physique Th\'eorique et Hautes Energies**} 
\centerline{Universit\'e de Paris XI, B\^atiment 210,
F-91405 Orsay Cedex, France}

\vskip 5 truemm
\centerline{\bf G. Velo}
\centerline{Dipartimento di Fisica}
\centerline{Universit\`a di Bologna and INFN, Sezione di Bologna, Italy}
\vskip 1 truecm
\begin{abstract}
\baselineskip=24 pt 
We study the theory of scattering in the energy space for the Hartree equation in space dimension
$n \geq 3$. Using the method of Morawetz and Strauss, we prove in par\-ti\-cu\-lar asymptotic
completeness for radial nonnegative nonincreasing potentials satisfying suitable regularity
properties at the origin and suitable decay properties at infinity. The results cover in particular
the case of the potential $|x|^{-\gamma}$ for $2 < \gamma < {\rm Min}(4,n)$.
   \end{abstract}
\baselineskip=12 pt
\vskip 3 truecm
\noi AMS Classification : Primary 35P25. Secondary 35B40, 35Q40, 81U99. \par \vskip 5 truemm
\noi Key words : Scattering theory, Hartree equation, wave operators, asymptotic completeness.
 \vskip 1 truecm

\noindent LPTHE Orsay 98-57 \par
\noindent September 1998 \par
\vfill
\noi * Work supported in part by NATO Collaborative Research Grant 972231. \par
\noi ** Laboratoire associ\'e
au Centre National de la Recherche Scientifique - URA D0063.

\newpage
\pagestyle{plain}
\baselineskip=24 pt

\mysection{Introduction.}
\hspace*{\parindent}
This paper is devoted to the theory of scattering for the Hartree equation
\beq
\label{1.1}
i\partial_t u = - (1/2)\Delta u + u\left ( V \star |u|^2 \right ) \quad .
\eeq  
\noi Here $u$ is a complex valued function defined in space time ${I\hskip-1truemm R}^{n+1}$,
$\Delta$ is the Laplacian in ${I\hskip-1truemm R}^n$, $V$ is a real valued even function defined
in ${I\hskip-1truemm R}^n$, hereafter called the potential, and $\star$ denotes the convolution in
${I\hskip-1truemm R}^n$. In particular we develop a complete theory of scattering for the equation
(\ref{1.1}) in the energy space, which turns out to be the Sobolev space $H^1$, under suitable
assumptions on $V$. \par

One of the basic problems addressed by scattering theory is that of classifying the asymptotic
behaviours in time of the global solutions of a given evolution equation by comparing
them with those of the solutions of a suitably chosen and simpler evolution equation. In the
special case where the given equation is a (nonlinear) perturbation of a linear dispersive
equation, the first obvious candidate as a comparison equation is the underlying linear
equation, hereafter called the free equation, and we restrict the discussion to that case. In the
case of the Hartree equation (\ref{1.1}), that equation is the free Schr\"odinger equation
\beq
\label{1.2}
i \partial_t u = - (1/2) \Delta u \quad .
\eeq  

Let $U(t)$ be the evolution group that solves the free equation. The comparison between the two
equations then gives rise to the following two basic questions. \par

{\bf (1) Existence of the wave operators.} For any solution $v_+(t) = U(t) \ u_+$ of the free
equation with initial data $v_+(0) = u_+$, called the asymptotic state, in a suitable Banach space
$Y$, one looks for a solution $u$ of the original equation which behaves asymptotically as $v_+$
when $t \to \infty$, typically in the sense that
\beq
\label{1.3}
\parallel u(t) - v_+(t);Y \parallel \ \to 0 \qquad \hbox{when} \ t \to + \infty
\eeq 
\noi or rather
\beq
\label{1.4}
\parallel U(-t) u(t) - u_+ ; Y \parallel \ \to 0 \qquad \hbox{when} \ t \to + \infty
\eeq
\noi which may be more appropriate if $U(\cdot )$ is not a bounded group in $Y$. If such a $u$ can
be constructed for any $u_+ \in Y$, one defines the wave operator $\Omega_+$ for positive time as
the map $u_+ \to u(0)$. The same problem arises at $t \to - \infty$, thereby
leading to the definition of the wave operator $\Omega_-$ for negative time. \par

{\bf (2) Asymptotic completeness.} Conversely, given a solution $u$ of the original equation, one
looks for asymptotic states $u_+$ and $u_-$ such that $v_{\pm}(t) = U(t) \ u_{\pm}$ behaves
asymptotically as $u(t)$ when $t \to \pm \infty$, typically in the sense that (\ref{1.3}) or
(\ref{1.4}) and their analogues for negative time hold. If that can be realized for any $u$ with
initial data $u(0)$ in $Y$ for some $u_{\pm} \in Y$, one says that asymptotic completeness holds in
$Y$. \par

Asymptotic completeness is a much harder problem than the existence of the wave operators, except
in the case of small data where it follows as an immediate by-product of the method generally used
to solve the latter problem. Asymptotic completeness for large data in the sense described above
requires strong assumptions on the nonlinear perturbation, in particular some repulsivity condition,
and proceeds through the derivation of a priori estimates for general solutions of the original
equation. It has been derived so far only for a small number of nonlinear evolution equations. \par

Scattering theory for nonlinear evolution equations started in the sixties under the impulse of
Segal \cite{24r} \cite{25r}, on the example of the nonlinear wave (NLW) equation
\beq
\label{1.5}
\sq u + f(u) = 0
\eeq
\noi and of the nonlinear Klein-Gordon (NLKG) equation
\beq
\label{1.6}
\left ( \sq + m^2 \right ) u + f(u) = 0
\eeq
\noi with $f(u)$ a nonlinear interaction term, a typical form of which is 
\beq
\label{1.7}
f(u) = |u|^{p-1} u
\eeq
\noi for some $p > 1$. Following early works where asymptotic completeness was ensured by including
an explicit integrable decay in time or a suitable decay in space in the nonlinearity (see
references quoted in \cite{26r}), major contributions were made by Strauss on the NLW equation
\cite{26r} and by Morawetz and Strauss on the NLKG equation \cite{22r}. In \cite{26r} a complete
theory of scattering is developed for the NLW equation (\ref{1.5}) in space dimension $n = 3$, in a
space of suitably regular and decaying functions, for a large class of nonlinearities including the
form (\ref{1.7}) under the natural assumption $3 \leq p < 5$ (more generally $4/(n-1) \leq p - 1 <
4/(n - 2)$ in space dimension $n$). Asymptotic completeness is proved there by exploiting the
approximate conservation law associated with the approximate conformal invariance of the equation. In
\cite{22r}, a complete theory of scattering is developed for the NLKG equation (\ref{1.6}) in space
dimension $n = 3$ with nonlinear interaction (\ref{1.7}) and $p = 3$. Asymptotic completeness 
is proved there by the use of the Morawetz inequality \cite{21r}. That inequality however is not a
very strong statement, since it asserts only the convergence of a space time integral which would
naively be expected to be barely divergent, and it required a tour de force in analysis to extract
therefrom the necessary a priori estimates. The method of \cite{22r} was then extended by Lin and
Strauss \cite{20r} to construct a complete theory of scattering for the nonlinear Schr\"odinger (NLS)
equation \beq
\label{1.8}
i \partial_t u = - (1/2) \Delta u + f(u)
\eeq
\noi in space dimension $n = 3$ with nonlinear interaction (\ref{1.7}) and $8/3 < p < 5$. Parallel
and subsequent developments included the construction of a complete theory of scattering in the
space $\Sigma = H^1 \cap {\cal F}H^1$, where $H^1$ is the usual Sobolev space and ${\cal F}$ the
Fourier transform, in arbitrary space dimension, both for the NLS equation (\ref{1.8}) \cite{9r}
\cite{17r} \cite{31r} and for the Hartree equation (\ref{1.1}) \cite{10r} \cite{18r} \cite{23r}.
Asymptotic completeness is proved there by the use of the approximate conservation law
associated with the approximate pseudo-conformal invariance of the NLS and Hartree equations, which
is the analogue for those equations of the conformal invariance of the NLW equation exploited in
\cite{26r}. The class of interactions thereby covered includes the form (\ref{1.7}) for the NLS
equation (\ref{1.8}) with  \beq
\label{1.9}
p_0 (n) < p < 1 + 4/(n - 2) \quad , 
\eeq
\noi where $p_0(n)$ is the positive root of the equation $np(p - 1) = 2(p + 1)$, and includes the
potential
\beq
\label{1.10}
V(x) = C \ |x|^{- \gamma}
\eeq
with $C > 0$ and $4/3 < \gamma < {\rm Min}(4,n)$ for the Hartree equation (\ref{1.1}). \par

Meanwhile the paper \cite{22r} inspired further developments. It turned out that the natural
function space of initial data and asymptotic states for the implementation of the method of
\cite{22r} is the energy space, and a complete theory of scattering in that space was
constructed for the NLKG equation (\ref{1.6}) in arbitrary dimension $ n \geq 3$, under assumptions
on $f$ which in the special case (\ref{1.7}) reduce to the natural condition 
\beq
\label{1.11}
1 + 4/n < p < 1 + 4/(n - 2)
\eeq
\noi \cite{5r} \cite{6r}, see also \cite{12r}. A complete theory of scattering in the energy space,
in that case the Sobolev space $H^1$, was then constructed for the NLS equation in dimension $n
\geq 3$ again by the use of a variant of the method of \cite{22r}, under assumptions on $f$ which
in the special case (\ref{1.7}) again reduce to the natural condition (\ref{1.11}) \cite{11r}
\cite{12r}. The construction was then somewhat simplified in \cite{7r} and in \cite{8r}. Finally,
the method of \cite{22r} was extended to construct a complete theory of scattering in the energy
space for the NLW equation (\ref{1.5}) under assumptions on $f$ which barely miss the special case
(\ref{1.7}) with $p = 1 + 4/(n - 2)$, the $H^1$ critical value \cite{13r}. However the proof of
asymptotic completeness by that method for that equation is now superseded by more direct
estimates which cover that critical case \cite{1r} \cite{2r} \cite{3r}. Expositions of the theory
at various stages of its development can be found in \cite{27r} \cite{28r} \cite{30r}. \par

In this paper, we develop a complete theory of scattering for the Hartree equation (\ref{1.1}) in
the energy space, which is again the Sobolev space $H^1$. The essential part of that theory is the
proof of asymptotic completeness, which is obtained by an adaptation of the method of \cite{22r}.
As a preliminary, we present briefly the theory of the Cauchy problem at finite times and the
construction of the wave operators in the energy space $H^1$. That part of the theory is a simple
variant of the corresponding theory for the NLS equation, which has reached a well developed
stage \cite{7r} \cite{8r} \cite{19r}. Consequently, although we give complete statements of the
results, we provide only brief sketches of the proofs in that part. The exposition follows
closely that in \cite{8r} for the NLS equation. \par

In all this paper, we restrict our attention to space dimension $n \geq 3$, because the method
of \cite{22r} applies only to that case. Most of the results on the Cauchy problem at finite
times and on the existence of the wave operators would extend to lower dimensions, where they
would however require modified statements. \par

The assumptions made on $V$ for the Cauchy problem at finite times and for the existence of the
wave operators are of a general character. In the typical situation where $V \in L^p$ for some
$p$, $1 \leq p \leq \infty$, they reduce to the condition $p > n/4$ for the Cauchy problem at
finite times and to $n/4 < p \leq n/2$ for the existence of the wave operators, thereby covering
the example (\ref{1.10}) for $\gamma < {\rm Min} (4, n)$ and $2 < \gamma < {\rm Min}(4, n)$
respectively. On the other hand stronger assumptions are required for the proof of asymptotic
completeness. In particular the potential $V$ should in addition be radial and suitably
repulsive (see Assumption (H3) in Section 4 below). Those requirements are satisfied in
particular by the special case (\ref{1.10}), and a complete theory can be developed for that case
for $2 < \gamma < {\rm Min} (4, n)$. \par

This paper is organized as follows. In Section 2 we treat the Cauchy problem at finite times for
the equation (\ref{1.1}). We prove local wellposedness in $H^1$ (Proposition 2.1), we derive
the conservation laws of the $L^2$ norm and of the energy (Proposition 2.2), and we prove
global wellposedness in $H^1$ (Proposition 2.3). In Section 3, we prove the existence of the
wave operators. We solve the local Cauchy problem in a neighbourhood of infinity in time
(Proposition 3.1), we prove the existence and some properties of asymptotic states for the
solutions thereby obtained (Proposition 3.2), and we conclude with the existence of the wave
operators (Proposition 3.3). Finally in Section 4, we prove the main result of this paper,
namely asymptotic completeness in $H^1$. We derive the Morawetz inequality for the Hartree
equation (\ref{1.1}) (Proposition 4.1), we extract therefrom an estimate of the solutions in
suitable norms (Proposition 4.2), and we exploit that estimate to prove asymptotic
completeness (Proposition 4.3). A more detailed description of the contents of Section 4 is
provided at the beginning of that section. \par

We conclude this introduction by giving some notation which will be used freely throughout this
paper. For any integer $n \geq 3$, we let $2^{\star} = 2n/(n-2)$. For any $r$, $1 \leq r \leq
\infty$, we denote by $\parallel \cdot \parallel_r$ the norm in $L^r \equiv L^r({I\hskip-1truemm
R}^n)$ and by $\bar{r}$ the conjugate exponent defined by $1/r + 1/\bar{r} = 1$, and we define
$\delta (r) = n/2 - n/r$. We denote by $<\cdot , \cdot >$ the scalar product in $L^2$ and by
$H_r^1$ the Sobolev space
$$H_r^1 = \left \{ u : \parallel u; H_r^1 \parallel \ = \ \parallel u \parallel_r \ + \
\parallel \nabla u \parallel_r \ < \infty \right \} \quad .$$
\noi For any interval $I$ of ${I\hskip-1truemm R}$, for any Banach space $X$, we denote by
${\cal C}(I, X)$ the space of continuous functions from $I$ to $X$ and for $1 \leq q \leq
\infty$, by $L^q (I, X)$ (resp. $L_{loc}^q(I, X)$) the space of measurable functions $u$ from
$I$ to $X$ such that $\parallel u(\cdot ) ; X\parallel \ \in L^q(I)$ (resp. $\in
L_{loc}^q(I)$). For any interval $I$ of ${I\hskip-1truemm R}$, we denote by $\bar{I}$ the
closure of $I$ in $\bar{{I\hskip-1truemm R}} = {I\hskip-1truemm R} \cup \{ \pm \infty \}$
equipped with the natural topology. Finally for any real numbers $a$ and $b$, we let $a \vee
b = {\rm Max} (a,b)$, $a \wedge b = {\rm Min} (a, b)$, $a_+ = a \vee 0$ and $a_- = (-a)_+$.

\mysection{The Cauchy problem at finite times.}
\hspace*{\parindent} 
In this section, we briefly recall the relevant results on the Cauchy problem with finite initial
time for the equation (\ref{1.1}). We refer to \cite{7r} \cite{10r} for more details. We rewrite
the equation (\ref{1.1}) as
\beq
i \partial_t u = - (1/2) \Delta u + f(u)
\label{2.1}
\eeq 
\noi where
\beq
\label{2.2}
f(u) = u \left ( V  \star  |u|^2 \right ) \quad .
\eeq
 \noi The Cauchy problem for the equation (\ref{2.1}) with initial data $u(0) = u_0$ at $t = 0$ is
rewritten in the form of the integral equation
\beq
\label{2.3}
u(t) = U(t) \ u_0 - i \int_0^t dt' \ U(t - t') \ f(u(t'))
\eeq
\noi where $U(t)$ is the unitary group

$$U(t) = \exp \left ( i t \Delta /2 \right ) \quad .$$
\noi It is well known that $U(t)$ satisfies the pointwise estimate
\beq
\label{2.4}
\parallel U(t) f \parallel_r \ \leq (2 \pi |t|)^{-\delta (r)} \ \parallel f\parallel_{\bar{r}}
\eeq
\noi where $2 \leq r \leq \infty$, $1/r + 1/\bar{r} = 1$ and $\delta (r) \equiv n/2 - n/r$. \par

Let now $I$ be an interval, possibly unbounded. We define the operators
\beq
\label{2.5}
(U \star f)(t) = \int_I dt' \ U(t-t') \ f(t') \quad ,
\eeq
\beq
\label{2.6}
(U \star_R f)(t) = \int_{I \cap \{t' \leq t\}} dt'\ U(t - t') \ f(t') \quad ,
\eeq
\noi where $f$ is defined in ${I\hskip-1truemm R}^n \times I$ and suitably regular, the
dependence on $I$ is omitted and the subscript $R$ stands for retarded. We introduce the following
definition \\

\noi {\bf Definition 2.1.} A pair of exponents $(q, r)$ is said to be admissible if
$$0 \leq 2/q = \delta (r) < 1 \quad .$$
\noi It is well known that $U(t)$ satisfies the following Strichartz estimates \cite{7r} \cite{8r}
\cite{19r}. \\

\noi {\bf Lemma 2.1.} {\it The following estimates hold~: \par
(1) For any admissible pair $(q,r)$}
\beq
\label{2.7}
\parallel U(t) \ u ; L^q({I\hskip-1truemm R}, L^r) \parallel \ \leq c_r \ \parallel u \parallel_2
\quad . \eeq 

{\it (2) For any admissible pairs $(q_i, r_i)$, $i = 1,2$, and for any interval $I \subset {I\hskip-1truemm
R}$}
\beq
\label{2.8}
\parallel U \star f;L^{q_1}(I, L^{r_1}) \parallel \ \leq c_{r_1} \ c_{r_2} \ \parallel
f;L^{\bar{q}_2}(I, L^{\bar{r}_2})\parallel \quad ,  \eeq
\beq
\label{2.9}
\parallel U \star_R f;L^{q_1}(I, L^{r_1}) \parallel \ \leq c_{r_1} \ c_{r_2} \ \parallel
f;L^{\bar{q}_2}(I, L^{\bar{r}_2})\parallel \quad .  \eeq
\vskip 5 truemm

Lemma 2.1 suggests that we study the Cauchy problem for the equation (\ref{2.3}) in spaces of the
following type. Let $I$ be an interval. We define
\beq
\label{2.10}
X(I) = \left \{ u:u \in {\cal C}(I,L^2) \ \hbox{and} \ u \in L^q(I, L^r) \ \hbox{for} \ 0 \leq 2/q
= \delta (r) < 1 \right \} \quad ,  \eeq
\beq
\label{2.11}
X^1(I) = \left \{ u : u \ \hbox{and} \ \nabla u \in X(I) \right \} \quad .
\eeq

For noncompact I, we define the spaces $X_{loc}(I)$ and $X^1_{loc}(I)$ in a similar way by
replacing $L^q$ by $L^q_{loc}$ in (\ref{2.10}). \par

The spaces $X(I)$ and $X^1(I)$ are not Banach spaces in a natural way because the interval for $r$ in
(\ref{2.10}) is semi-open. We define Banach spaces by restricting the interval for $r$ by $0 \leq
\delta (r) \leq 2/q_0 = \delta (r_0) \equiv \delta_0 < 1$, namely
\bea
\label{2.12}
X_{r_0}(I) &=& \left \{ u:u \in {\cal C} (I, L^2) \ \hbox{and} \ u \in L^q(I, L^r) \ \hbox{for} \ 0
\leq 2/q = \delta (r) \leq \delta_0 \right \} \nn \\
&=& ({\cal C} \cap L^{\infty}) (I, L^2) \cap L^{q_0} (I, L^{r_0}) \quad , \\    
X^1_{r_0}(I) &=& \left \{ u : u \ \hbox{and} \ \nabla u \in X_{r_0}(I) \right \} \quad .
\label{2.13}
\eea
\noi The spaces $X_{r_0,loc}(I)$ and $X_{r_0,loc}^1(I)$ are defined in a natural way. \par

In all this paper, we assume that the potential $V$ satisfies the following assumption. \\

\noi (H1) $V$ is a real even function and $V \in L^{p_1} + L^{p_2}$ for some $p_1$, $p_2$ satisfying
\beq
\label{2.14}
1 \vee (n/4) \leq p_2 \leq p_1 \leq \infty \quad .
\eeq

We can now state the main result on the local Cauchy problem for the equation (\ref{1.1}) with
$H^1$ initial data. \\

\noi {\bf Proposition 2.1.} {\it Let $V$ satisfy (H1) and define $r_0$ by}
\beq
\label{2.15}
\delta_0 \equiv \delta (r_0) = (n/4p_2 - 1/2)_+ \quad (\leq 1/2) \quad .
\eeq
\noi {\it Let $u_0 \in H^1$. Then \par
(1) There exists a maximal interval $(-T_-, T_+)$ with $T_{\pm} > 0$ such that the equation
(\ref{2.3}) has a unique solution $u \in X_{r_0,loc}^1((-T_-, T_+))$. The solution $u$ actually
belongs to $X_{loc}^1((-T_-, T_+))$. \par

(2) For any interval $I$ containing 0, the equation (\ref{2.3}) has at most one solution in
$X_{r_0}^1(I)$. \par

(3) For $- T_- < T_1 \leq T_2 < T_+$, the map $u_0 \to u$ is continuous from $H^1$ to $X^1([T_1,
T_2])$. \par

(4) Let in addition $p_2 > n/4$. Then if $T_+ < \infty$ (resp. $T_- < \infty$), $\parallel
u(t);H^1\parallel \to \infty$ when $t$ increases to $T_+$ (resp. decreases to $-T_-$).} \\

\noi {\bf Sketch of proof.} The proof proceeds by standard arguments. The main technical point
consists in proving that the operator defined by the RHS of (\ref{2.3}) is a contraction in
$X_{r_0}^1(I)$ on suitable bounded sets of $X_{r_0}^1(I)$ for $I = [-T, T]$ and $T$ sufficiently
small. The basic estimates follow from Lemma 2.1 and from
\beq
\label{2.16}
\parallel f(u) ; L^{\bar{q}}(I, H_{\bar{r}}^1)\parallel \ \leq C \parallel V \parallel_p \
\parallel u;L^q(I,H_r^1)\parallel \ \parallel u;L^k(I,L^s)\parallel^2 \ T^{\theta} \eeq 

\noi where we have assumed for simplicity that $V \in L^p$, where $(q, r)$ is an admissible pair,
and where the exponents satisfy
\beq
\label{2.17}
n/p = 2 \delta (r) + 2 \delta (s)
\eeq
\beq
\label{2.18}
2/q + 2/k = 1 - \theta \quad .
\eeq

The estimate (\ref{2.16}) is obtained by applying the H\"older and Young inequalities in space
followed by the H\"older inequality in time, which requires $0 \leq \theta \leq 1$. A similar
estimate holds for the difference of two solutions. For general $V$ satisfying (H1), the
contribution of the components in $L^{p_1}$ and $L^{p_2}$ are treated separately, with the
exponents $(q,r,k,s)$ possibly depending on $p$. \par

If $n/p \leq 2$, one can choose $r = 2$, $\delta (s) = n/2p \leq 1$, $k = q = \infty$, so that
$\theta = 1$ and one can take $r_0 = 2$. \par

If $n/p \geq 2$, one can take $\delta (s) = \delta (r) + 1$ and $k = q$, so that $4/q = 1 - \theta$
and $n/p = 4\delta (r) + 2 = 4 -2\theta$, which yields $\theta \geq 0$ for $n/p \leq 4$, and allows
for $\delta_0 = n/4p-1/2$. \par

The $H^1$-critical case $p = n/4$ yields $\theta = 0$ and requires a slightly more refined
treatment than the subcritical case $p > n/4$. \par \nobreak
\hfill $\sq$ \\

It is well known that the Hartree equation (\ref{1.1}) formally satisfies the conservation of the
$L^2$ norm and of the energy
\beq
\label{2.19}
E(u) = {1 \over 2} \parallel \nabla u \parallel_2^2 \ + {1 \over 2} \int dx \ dy \ |u(x)|^2 \ V(x -
y)\ |u(y)|^2 \quad .  \eeq
\noi Actually it turns out that the $X^1$ regularity of the solutions constructed in Proposition
2.1 is sufficient to ensure those conservation laws. \\

\noi {\bf Proposition 2.2.} {\it Let $V$ satisfy (H1) and define $r_0$ by (\ref{2.15}). Let $I$ be
an interval and let $u \in X_{r_0}^1(I)$ be a solution of the equation (\ref{1.1}). Then $u$
satisfies $\parallel u(t_1)\parallel_2 \ = \ \parallel u(t_2)\parallel_2$ and $E(u(t_1)) = E(u(t_2))$
for all $t_1$, $t_2 \in I$.} \\

\noi {\bf Sketch of proof.} We consider only the more difficult case of the energy and we follow
the proof of the corresponding result for the nonlinear Schr\"odinger (NLS) equation given in
\cite{8r}. Let $\varphi \in {\cal C}_0^{\infty}$ be a smooth approximation of the Dirac
distribution $\delta$ in ${I\hskip-1truemm R}^n$. By an elementary computation which is allowed by
the available regularity, one obtains
\bea
\label{2.20}
&&E\left ( (\varphi \star u)(t_2) \right ) - E \left ( (\varphi \star u)(t_1)\right ) = -
{\rm Im} \int_{t_1}^{t_2} dt \Big \{ <\varphi \star \nabla u, \nabla f(\varphi \star u) - \varphi \star \nabla f(u)> 
\nn \\
&&+ 2<\varphi \star f(u), f(\varphi \star u) - \varphi 
\star f(u)> \Big \} (t) \quad . \eea
\noi One then lets $\varphi$ tend to $\delta$, using the fact that convolution with $\varphi$ tends
strongly to the unit operator in $L^r$ for $1 \leq r < \infty$. The LHS of (\ref{2.20}) tends to
$E(u(t_2)) - E(u(t_1))$ and the RHS is shown to tend to zero by the Lebesgue dominated convergence
theorem applied to the time integration. For that purpose one needs an estimate of the integrand
which is uniform in $\varphi$ and integrable in time. That estimate essentially boils down to
\beq
\label{2.21}
|<\nabla u, \nabla f>| \ \leq C \parallel V \parallel_p \ \parallel \nabla u \parallel_r^2 \
\parallel u \parallel_s^2 \eeq
\noi with $r$ and $s$ satisfying (\ref{2.17}) and to
\beq
\label{2.22}
\parallel f(u) \parallel_2^2 \ \leq C \parallel V \parallel_p^2 \ \parallel u \parallel_s^6
\eeq
\noi with $\delta (s) = n/3p$, for a potential $V \in L^p$. \par

For (\ref{2.21}), we choose the same values of $r$, $s$ as in the proof of Proposition 2.1, so that
the RHS of (\ref{2.21}) belongs to $L^{\infty}$ in time for $n/p \leq 2$ and to $L^{q/4}$ with $q
\geq 4$ for $n/p \geq 2$. On the other hand, the RHS of (\ref{2.22}) belongs to $L^{\infty}$ in
time for $n/p \leq 3$ and to $L^{q/6}$ with $2/q = \delta (s) - 1 = n/3p - 1 \leq 1/3$ for $3 \leq
n/p \leq 4$. \par \nobreak
\hfill $\sq$ \\

We now turn to the global Cauchy problem for the equation (\ref{1.1}). For that purpose we need to
ensure that the conservation of the $L^2$ norm and of the energy provides an a priori estimate of
the $H^1$ norm of the solution. This is the case if the potential $V$ satisfies the following
assumption \\

\noi (H2) $V_- \equiv V \wedge 0 \in L^{n/2} + L^{\infty}$. \\

In fact, it follows from (H2) by the H\"older, Young and Sobolev inequalities that for any
$\varepsilon > 0$, there exists $C(\varepsilon )$ such that
\beq
\label{2.23}
\int dx \ dy \ |u(x)|^2 \ V_-(x - y) |u(y)|^2 \leq \varepsilon \parallel u \parallel_2^2 \
\parallel \nabla u \parallel_2^2 + C(\varepsilon ) \parallel u \parallel_2^4
 \eeq
\noi and therefore
\beq
\label{2.24}
\parallel \nabla u \parallel_2^2 \ \leq 4 E(u) + 2C\left ( (2 \parallel u \parallel_2^2 )^{-1} \right
) \ \parallel u \parallel_2^4 \quad . \eeq
\noi We can now state the main result on the global Cauchy problem for the equation (1.1). \\

\noi {\bf Proposition 2.3.} {\it Let $V$ satisfy (H1) with $p_2 > n/4$ and (H2). Let $u_0 \in H^1$
and let $u$ be the solution of the equation (\ref{2.3}) constructed in Proposition 2.1. Then $T_+ =
T_- = \infty$ and $u \in X_{loc}^1({I\hskip-1truemm R}) \cap L^{\infty} ({I\hskip-1truemm
R},H^1)$.} \\

The proof is standard. Note however that the result is stated only for the $H^1$ subcritical case
$p_2 > n/4$.

\mysection{Scattering Theory I. Existence of the wave operators.}
\hspace*{\parindent} 
In this section we begin the study of the theory of scattering for the Hartree equation (\ref{1.1})
by addressing the first question raised in the introduction, namely that of the existence of the
wave operators. We restrict our attention to positive time. We consider an asymptotic state $u_+
\in H^1$ and we look for a solution $u$ of the equation (\ref{1.1}) which is asymptotic to the
solution $v(t) = U(t)u_+$ of the free equation. For that purpose, we introduce the solution
$u_{t_0}$ of the equation (\ref{1.1}) satisfying the initial condition $u_{t_0}(t_0) = v(t_0)
\equiv U(t_0)u_+$. We then let $t_0$ tend to $\infty$. In favourable circumstances, we expect
$u_{t_0}$ to converge to a solution $u$ of the equation (\ref{1.1}) which is asymptotic to $v(t)$.
The previous procedure is easily formulated in terms of integral equations. The Cauchy problem with
initial data $u(t_0)$ at time $t_0$ is equivalent to the equation
\beq
\label{3.1}
u(t) = U(t - t_0) \ u(t_0) - i \int_{t_0}^t dt' \ U(t - t') \ f(u(t')) \quad .
\eeq
\noi The solution $u_{t_0}$ with initial data $U(t_0)u_+$ at time $t_0$ should therefore be a
solution of the equation
\beq
\label{3.2}
u(t) = U(t) u_+ - i \int_{t_0}^t dt' \ U(t - t') \ f(u(t')) \quad .
\eeq
\noi The limiting solution $u$ is then expected to satisfy the equation 
\beq
\label{3.3}
u(t) = U(t) \ u_+ + i \int_t^{\infty} dt' \ U(t - t') \ f(u(t')) \quad .
\eeq
\noi The problem of existence of the wave operators is therefore the Cauchy problem with infinite
initial time. We solve that problem in two steps. We first solve it locally in a neighborhood of
infinity by a contraction method. We then extend the solutions thereby obtained to all times by
using the available results on the Cauchy problem at finite times. In order to solve the local
Cauchy problem at infinity, we need to use function spaces including some time decay in their
definition, so that at the very least the integral in (\ref{3.3}) converges at infinity.
Furthermore the free solution $U(t)u_+$ should belong to those spaces. In view of Lemma 2.1,
natural candidates are the spaces $X_{r_0}^1(I)$ for some $I = [T, \infty )$, where the time decay is
expressed by the $L^q$ integrability at infinity, and we shall therefore study that problem in
those spaces. We shall also need the fact that the time decay of $u$ implies sufficient time
decay of $f(u)$. This will show up through additional assumptions on $V$ in the form of an upper
bound on $p_1$, namely $p_1 \leq n/2$. \par

For future reference we state additional time integrability properties of functions in
$X_{r_0}^1(\cdot )$ which are not immediately apparent on the definition. \\

\noi {\bf Lemma 3.1.} {\it Let $I$ be an interval, possibly unbounded. Then}
\beq
\label{3.4}
\parallel u;L^{q_0}(I, L^r) \parallel \ \leq C \parallel u;X_{r_0}^1(I) \parallel
\eeq   
\noi {\it for $\delta_0 \leq \delta (r) \leq \delta_0 + 1$, where $C$ is independent of $I$.} \\

\noi {\bf Proof.} The result follows from the Sobolev inequality 
$$\parallel u \parallel_r \ \leq C \parallel u \parallel_{r_0}^{1 - \sigma} \ \parallel \nabla u
\parallel_{r_0}^{\sigma}$$
\noi with $\sigma = \delta (r) - \delta (r_0)$ and from the definition of $X_{r_0}^1$. \par \nobreak
\hfill $\sq$ \\

We shall use freely the notation $\widetilde{u}(t) = U(-t) u(t)$ for $u$ a suitably regular
function of space time. We also recall the notation $\overline{{I\hskip-1truemm R}}$ for ${I\hskip-1truemm
R} \cup \{ \pm \infty \}$ and $\bar{I}$ for the closure of an interval $I$ in $\overline{{I\hskip-1truemm
R}}$ equipped with the obvious topology. \par

We can now state the main result on the local Cauchy problem in a neighborhood of infinity. \\

\noi {\bf Proposition 3.1.} {\it Let $r_0 = 2n/(n - 1)$, so that $\delta_0 = 1/2$. Let $V$
satisfy (H1) with $p_1 \leq n/2$. Let $u_+ \in H^1$. Then\par
(1) There exists $T < \infty$ such that for any $t_0 \in \bar{I}$ where $I \in [T, \infty )$, the
equation (\ref{3.2}) has a unique solution $u$ in $X_{r_0}^1(I)$. The solution $u$ actually belongs
to $X^1(I)$. \par
(2) For any $T' > T$, the solution $u$ is strongly continuous from $u_+ \in H^1$ and $t_0 \in
\bar{I}'$ to $X^1(I')$, where $I' = [T', \infty )$.} \\

\noi {\bf Sketch of proof.} The proof proceeds by a contraction argument in $X_{r_0}^1(I)$. The
main technical point consists in proving that the operator defined by the RHS of (\ref{3.2}) is a
contraction in $X_{r_0}^1(I)$ on suitable bounded sets of $X_{r_0}^1(I)$ for $T$ sufficiently large.
The basic estimate is again (\ref{2.16}) supplemented by (\ref{2.17}) (\ref{2.18}), now however
with $\theta = 0$. The fact that we use spaces where the time decay appears in the form of an
$L^q$ integrability condition in time forces the condition $\theta = 0$, so that we are in a
critical situation, as was the case for the local Cauchy problem at finite times in the $H^1$
critical case $p_2 = n/4$. We choose the exponents in (\ref{2.16}) as follows. We take $q = k =
q_0 = 4$, $\delta (r) = \delta_0 = 1/2$ and $1 + 2 \delta (s) = n/p$ so that the last norm in
(\ref{2.16}) is controlled by the $X_{r_0}^1$ norm for $1/2 \leq \delta (s) \leq 3/2$, namely $2
\leq n/p \leq 4$. The smallness condition which ensures the contraction takes the form 
\beq
\label{3.5}
\parallel U(t) \ u_+;L^{q_0}(I, H_{r_0}^1) \parallel \ \leq R_0
\eeq  

\noi for some absolute constant $R_0$. In particular the time $T$ of local resolution cannot be
expressed in terms of the $H^1$ norm of $u_+$ alone, as is typical of a critical situation. \par

The continuity in $t_0$ up to and including infinity follows from an additional application of
the same estimates. \par \nobreak
\hfill $\sq$ \\

An immediate consequence of the estimates in the proof of Proposition 3.1 is the existence of
asymptotic states for solutions of the equation (\ref{1.1}) in $X_{r_0}^1([T, \infty ))$ for some
$T$. Furthermore the conservation laws of the $L^2$ norm and of the energy are easily extended to
infinite time for such solutions. \\

\noi {\bf Proposition 3.2.} {\it Let $r_0 = 2n/(n - 1)$. Let $V$ satisfy (H1) with $p_1 \leq
n/2$. Let $T \in {I\hskip-1truemm R}$, $I = [T, \infty )$ and let $u \in X_{r_0}^1(I)$ be a
solution of the equation (\ref{1.1}). Then \par
(1) $\widetilde{u} \in {\cal C}(\bar{I} , H^1)$. In particular the following limit exists}
\beq
\label{3.6}
\widetilde{u}(\infty ) = \lim_{t \to \infty} \widetilde{u}(t)
\eeq   
\noi {\it as a strong limit in $H^1$. \par
(2) $u$ satisfies the equation (\ref{3.3}) with $u_+ = \widetilde{u}(\infty )$. \par
(3) $u$ satisfies the conservation laws}
\beq
\label{3.7}
\parallel \widetilde{u} (\infty ) \parallel_2 \ = \ \parallel u \parallel_2 \qquad , \quad
(1/2)\parallel \nabla \widetilde{u}(\infty ) \parallel_2^2 \ = E(u) \quad . \eeq

\vskip 5 truemm

\noi {\bf Sketch of proof.} {\bf Part (1).} We estimate for $T \leq t_1 \leq t_2$
\bea
\label{3.8}
&&\parallel \widetilde{u}(t_2) - \widetilde{u}(t_1) ; H^1 \parallel \ = \ \parallel \int_{t_1}^{t_2}
dt \ U(t_2 - t) \ f(u(t));H^1 \parallel \nn \\ 
&&\leq \ \parallel U \star f;X_{r_0}^1([t_1, t_2])\parallel \ \leq C \parallel
f(u);L^{\bar{q}_0}([t_1, t_2];H_{\bar{r}_0}^1)\parallel \eea
\noi and we estimate the last norm as in the proof of Proposition 3.1, namely by (\ref{2.16}) with
$\theta = 0$ and with the same choice of exponents. \par

\noi {\bf Part (2)} follows from Part (1) and from Proposition 3.1, especially part (2).\par

\noi {\bf Part (3).} From the conservation laws at finite time and from Part (1), it follows that
the following limits exist
\bea
\label{3.9}
\parallel \widetilde{u}(\infty )\parallel_2 &=& \lim_{t \to \infty} \parallel u(t) \parallel_2 \ =
\ \parallel u \parallel_2 \quad , \nn \\
\lim_{t \to \infty} P(u(t)) &=& E(u) - (1/2) \lim_{t \to \infty} \parallel \nabla \widetilde{u} (t)
\parallel_2^2 \nn \\ 
&=& E(u) - (1/2) \parallel \nabla \widetilde{u}(\infty ) \parallel_2^2 
\eea
\noi where
\beq
\label{3.10}
P(u) = {1 \over 2} \int dx \ dy |u(x)|^2 \ V(x - y) \ |u(y)|^2 \quad .
\eeq
\noi On the other hand
\beq
\label{3.11}
|P(u)| \leq C\parallel V \parallel_p \ \parallel u \parallel_r^4 \ \in L^{q_0/4} = L^1
\eeq
\noi by the H\"older and Young inequalities and by Lemma 3.1 with $\delta_0 = 1/2 \leq \delta (r) =
n/4p \leq 1$. \par

It then follows from (\ref{3.11}) that the limit in (\ref{3.9}) is zero. \par \nobreak
\hfill $\sq$ \\

The existence and the properties of the wave operators now follow from the previous local results
at infinity and from the global results of Section 2. \\

\noi {\bf Proposition 3.3.} {\it Let $r_0 = 2n/(n-1)$. Let $V$ satisfy (H1) and (H2) with $p_1
\leq n/2$ and $p_2 > n/4$. Then \par
(1) For any $u_+ \in H^1$, the equation (\ref{3.3}) has a unique solution $u$ in $X_{r_0, loc}^1({I\hskip-1truemm
R})$ with restriction in $X_{r_0}^1({I\hskip-1truemm R}^+)$. In addition, $u \in
X_{loc}^1({I\hskip-1truemm R})$ with restriction in $X^1({I\hskip-1truemm R}^+)$, and $\widetilde{u}
\in {\cal C}({I\hskip-1truemm R} \cup \{ + \infty \}, H^1)$. Furthermore $u$ satisfies the
conservation laws} $$\parallel u(t) \parallel_2 \ = \ \parallel u_+\parallel_2 \qquad , \quad
E(u(t)) = (1/2) \parallel \nabla u_+ \parallel_2^2$$
\noi {\it for all $t \in {I\hskip-1truemm R}$.\par
(2) The wave operator $\Omega_+ : u_+ \to u(0)$ is well defined in $H^1$, and is continuous and
bounded in the $H^1$ norm.} \\

\noi {\bf Sketch of proof.} Part (1) follows immediately from Propositions 2.2, 2.3, 3.1 and 3.2.
In Part (2), boundedness of $\Omega_+$ follows from the conservation laws, while continuity follows
from the corresponding statements in Propositions 2.1 and 3.1. \par \nobreak
\hfill $\sq$ \\

The solutions of the equation (\ref{1.1}) constructed in Proposition 3.3, part 1 are dispersive at
$+ \infty$, but no claim is made at this stage on their behaviour at $- \infty$. Dispersiveness at
$- \infty$ would be a consequence of asymptotic completeness, which will be studied only in Section
4.

\mysection{Scattering theory II. Asymptotic completeness.}
\hspace*{\parindent}
In this section, we continue the study of the theory of scattering for the Hartree equation
(\ref{1.1}) by addressing the second question raised in the introduction, namely that of asymptotic
completeness. In particular we prove the main result of this paper, namely the fact that asympotic
completeness holds in the energy space $H^1$ for radial and suitably repulsive potentials (see
Assumption (H3) below). In view of the results of Section 3, especially Proposition 3.2, it will
turn out that the crux of the argument consists in showing that the global solutions of the
equation (\ref{1.1}) in $X_{loc}^1({I\hskip-1truemm R})$ constructed in Proposition 2.3 actually
belong to $X^1({I\hskip-1truemm R})$, namely exhibit the time decay properties contained in the
definition of that space. The proof uses the method of Morawetz and Strauss \cite{22r} and relies
on two estimates. The first one is an elementary propagation estimate which for the Hartree (as
well as for the NLS) equation replaces the finiteness of the propagation speed for the NLKG
equation (see Lemma 4.2). The second estimate follows from the Morawetz inequality, which is
closely related to the approximate dilation invariance of the equation (see Proposition 4.1).
Space time is split into an internal and an external region where $|x|$ is small or large
respectively as compared with $|t|$. For radial repulsive potentials according to the assumption
(H3) below, the Morawetz inequality implies an a priori estimate for a suitable norm of the
internal part of $u$ (see Proposition 4.2 and Lemma 4.4). One uses that estimate in the internal
region and the propagation estimate in the external region. Plugging those estimates into the
integral equation for the solution $u$, one proves successively that a suitable norm of $u$ is
small in large intervals (see Lemma 4.5) and tends to zero at infinity (see Lemma 4.6) and that
$u$ belongs to $X^1({I\hskip-1truemm R})$ (see Proposition 4.3). \par

We continue to assume $n \geq 3$ as in the rest of this paper. We restrict our attention to
positive times. We first state an elementary property of $H^1$ solutions of the free
Schr\"odinger equation. We recall that $2^{\star} \equiv 2n/(n-2)$. \\

\noi {\bf Lemma 4.1.} {\it Let $u_0 \in H^1$. Let $2 < r \leq 2^{\star}$. Then $U(t) u_0$ tends
to zero in $L^r$ norm when $|t| \to \infty$.} \\

\noi {\bf Proof.} We approximate $u_0$ in $H^1$ norm by $u'_0 \in L^{\bar{r}} \cap H^1$. By
(\ref{2.4}), Sobolev inequalities and the unitarity of $U(t)$ in $H^1$, we estimate 
\bea
\label{4.1}
\parallel U(t) \ u_0 \parallel_r &\leq& \parallel U(t) \ u'_0 \parallel_r \  + \ C \parallel u_0 -
u'_0 \parallel_2^{1 - \delta (r)} \ \parallel \nabla (u_0 - u'_0) \parallel_2^{\delta (r)} \nn \\
&\leq& (2 \pi |t|)^{- \delta (r)} \ \parallel u'_0\parallel_{\bar{r}} \ + \ C \parallel u_0 - u'_0 ;
H^1 \parallel \eea
\noi from which the result follows. \par \nobreak
\hfill $\sq$ \\

We next state the propagation property of finite energy solutions of the equation (\ref{1.1})
mentioned previously. For any function $u$ of space time and for $t \geq 1$, we define 
\beq
\label{4.2}
u_{>\atop<}(t,x) = u(t, x) \ \chi (|x| \mathrel{\mathop >_{\displaystyle{<}}} \ t \ {\rm Log} \ t)
\eeq   
\noi so that $u = u_> + u_<$. This decomposition corresponds to the splitting of space time
mentioned previously. There is nothing magic about the function $t \ {\rm Log} \ t$. It is chosen so
as to tend to infinity faster than $t$ and to ensure the divergence of the
integral $\int_{\cdot}^{\infty} dt (t \ {\rm Log} t)^{-1}$.\\

\noi {\bf Lemma 4.2.} {\it Let $V$ satisfy (H1), let $u \in ({\cal C} \cap L^{\infty})
({I\hskip-1truemm R}, H^1)$ be a solution of the equation (\ref{1.1}) and let $u_0 = u(0)$. Then \par
(1) For any $R > 0$ and any $t \in {I\hskip-1truemm R}$, $u$ satisfies the estimate}
\beq
\label{4.3}
\int_{|x| \geq R} dx |u(t, x)|^2 \leq \int dx (1 \wedge R^{-1}|x|)|u_0(x)|^2 + R^{-1} |t|
\parallel u \parallel_2 \ \parallel \nabla u; L^{\infty}({I\hskip-1truemm R}, L^2)\parallel
\quad . \eeq

{\it (2) For any $r$ with $2 \leq r < 2^{\star}$, $\parallel u_>(t)\parallel_r$ tends to zero when
$t \to \infty$.} \\

\noi {\bf Proof.} {\bf Part (1).} We give only the formal computation, which is easily justified at
the available level of regularity. For $h$ a suitably smooth real function, we compute 
\beq
\label{4.4}
\partial_t <u,hu> \ = {\rm Im} \ <u, \nabla h \cdot \nabla u>
\eeq
\noi and therefore
\beq
\label{4.5}
<u(t), h \ u(t)> \ \leq \ <u_0, h \ u_0> + \parallel \nabla h \parallel_{\infty} \ \parallel u
\parallel_2 \int_0^t dt' \parallel \nabla u(t')\parallel_2 \eeq
\noi from which (\ref{4.3}) follows by taking $h(x) = 1 \wedge R^{-1} |x|$. \par

\noi {\bf Part (2)} for $r = 2$ follows from (\ref{4.3}) with $R = t \ {\rm Log} \ t$ and from the
Lebesgue dominated convergence theorem applied to the term containing $u_0$. The result for general
$r$ follows by interpolation between that for $r = 2$ and uniform boundedness in $L^{2^{\star}}$.
\par \nobreak
\hfill $\sq$ \\

The second main ingredient of the proof is the Morawetz inequality which for the Hartree equation
can be written as follows. \\

\noi {\bf Proposition 4.1.} {\it Let $V$ satisfy (H1) and let $u \in X_{loc}^1({I\hskip-1truemm R})
\cap L^{\infty}({I\hskip-1truemm R} , H^1)$ be a solution of the equation (\ref{1.1}). Then for any
$t_1$, $t_2 \in {I\hskip-1truemm R}$, $t_1 \leq t_2$, the following inequality holds}
\beq
\label{4.6}
- \int_{t_1}^{t_2} dt \int dx \ \rho (x) \ \widehat{x}\cdot ( V \star \nabla \rho ) \leq 2
\parallel u \parallel_2 \ \parallel \nabla u ; L^{\infty} ({I\hskip-1truemm R}, L^2) \parallel \eeq 
\noi {\it where $\widehat{x} = |x|^{-1}x$ and $\rho = |u|^2$.} \\

\noi {\bf Proof.} There are several proofs of the Morawetz inequality for various equations in the
literature. Here we follow the version given in \cite{8r} \cite{11r}. In order to derive the result
at the available level of regularity, we introduce the same regularization as in the proof of
energy conservation in Proposition 2.2. We denote $\varphi \star u = u_{\varphi}$. Let $h$ be
a ${\cal C}^4$ function of space with bounded derivatives up to order 4. We compute 
\beq
\label{4.7}
\partial_t \ {\rm Im} <u_{\varphi} , (\nabla h) \cdot \nabla u_{\varphi}> \ = {\rm Re} \ <(\Delta h)
u_{\varphi} + 2(\nabla h) \cdot \nabla u_{\varphi}, - (1/2) \Delta u_{\varphi} + \varphi \star 
f(u) > \quad . \eeq
\noi The kinetic part of the RHS (namely the terms not containing $f$) is treated exactly as for
the NLS equation in \cite{11r}. The term containing $f$ is rewritten by using the fact that  
\beq
\label{4.8}
{\rm Re} \ <(\Delta h)u_{\varphi} + 2(\nabla h)\cdot \nabla u_{\varphi}, f(u_{\varphi})> \ = - \int
dx \ \rho_{\varphi}(\nabla h) \cdot (V \star \nabla \rho_{\varphi}) \eeq
\noi where $\rho_{\varphi} = |u_{\varphi}|^2$, and we obtain
\bea
\label{4.9}
&&\partial_t \ {\rm Im} <u_{\varphi},(\nabla h) \cdot \nabla u_{\varphi}> \ = \ <\nabla_i \
u_{\varphi}, (\nabla_{ij}^2h) \nabla_j u_{\varphi}> - (1/4) <u_{\varphi}, (\Delta^2h)u_{\varphi}> \nn
\\ &&- \int dx \ \rho_{\varphi} (\nabla h) \cdot (V \star \nabla \rho_{\varphi}) + \ {\rm
Re} \ <(\Delta h)u_{\varphi} + 2(\nabla h) \cdot \nabla u_{\varphi}, \varphi \star f(u) -
f(u_{\varphi})> \eea
\noi where $\nabla_{ij}^2h$ is the matrix of second derivatives of $h$ and summation over the dummy
indices $i$ and $j$ is understood. \par

We next let $\varphi$ tend to $\delta$. We integrate (\ref{4.9}) in time in the interval $[t_1,
t_2]$ and we take the limit of the time integral of the RHS by using the Lebesgue dominated
convergence theorem in the time variable. For that purpose, we need an estimate of the integrand
which is uniform in $\varphi$ and integrable in time. Such an estimate is obvious for the kinetic
terms. For the terms containing $f$, it boils down to
\beq
\label{4.10}
|<\nabla u \cdot \nabla h, f(u)>| \leq C \parallel \nabla h \parallel_{\infty} \ \parallel V
\parallel_p \ \parallel \nabla u \parallel_2 \ \parallel u \parallel_s^3 \eeq
\noi with $\delta (s) = n/3p \leq 4/3$, so that the RHS of (\ref{4.10}) belongs to $L_{loc}^2$ in
time. \par

All terms then tend to the obvious limits and we obtain
\bea
\label{4.11}
\left . {\rm Im} \ <u, (\nabla h) \cdot \nabla u> \right |_{t_1}^{t_2} = \int_{t_1}^{t_2} dt \Big \{
< \nabla_i u, (\nabla_{ij}^2h) \nabla_j u > \nn \\
- (1/4) <u, (\Delta^2h)u> - \int dx \ \rho (\nabla h) \cdot (V \star \nabla \rho) \Big \} (t)
\quad .
 \eea

We next take $h = (x^2 + \sigma^2)^{1/2}$ for some $\sigma > 0$ and we compute 
\begin{eqnarray*}
&&\nabla h = h^{-1} x \\
&&\nabla_{ij}^2h = h^{-1} \left ( \delta_{ij} - h^{-2} \ x_i \ x_j \right ) \\
&&\Delta^2h = - (n - 1) (n - 3) h^{-3} - 6(n - 3) \sigma^2 h^{-5} - 15 \sigma^4 h^{-7} 
\end{eqnarray*}
\noi so that $\nabla_{ij}^2 h$ is a positive matrix and $\Delta^2 h$ is negative. We then obtain an
inequality by dropping the kinetic terms in the RHS of (\ref{4.11}). Taking in addition the
harmless limit $\sigma \to 0$ by the Lebesgue dominated convergence theorem, we obtain
\beq
\label{4.12}
\left . - \int_{t_1}^{t_2} dt \int dx \ \rho (x) \widehat{x} \cdot (V \star \nabla \rho ) \leq
{\rm Im} <u, \widehat{x} \cdot \nabla u > \right |_{t_1}^{t_2} \eeq
\noi from which (\ref{4.6}) follows immediately. \par \nobreak
\hfill $\sq$ \\

For sufficiently regular $V$, for instance for $V \in {\cal C}^1$ with compact support, the
integrand of the time integral in (\ref{4.12}) can be rewritten as
\beq
\label{4.13}
- \int dx \ \rho (x) \widehat{x} \cdot (V \star \nabla \rho ) = - (1/2) \int dx \ dy
(\widehat{x} - \widehat{y}) \cdot \nabla V(x - y) \rho (x) \ \rho (y) \eeq
\noi which is suggestive of algebraic manipulations to be made below\footnote{The contribution of
the Hartree interaction to the Morawetz inequality given in \protect{\cite{7r}} (Theorem 7.4.4) is
incorrect. The error appears in (7.4.13) through the application of the radial derivative to a
convolution product.}. \par

The assumptions on $V$ made so far are not stronger than those made in Section 2. In particular the
assumptions on $u$ made in Lemma 4.2 and Proposition 4.1 are ensured by the assumptions of
Proposition 2.3, namely (H1) (H2) and $p_2 > n/4$, for any initial data in $H^1$. In order to
proceed further, we need to exploit the fact that the LHS of (\ref{4.6}) controls some suitable
norm of $u$ and for that purpose we need a repulsivity condition on $V$. That condition takes the
following form. \\

\noi (H3) $V$ is radial and nonincreasing, namely $V(x) = v(|x|)$ where $v$ is nonincreasing in 
${I\hskip-1truemm R}^+$. Furthermore, for some $\alpha \geq 2$, $v$ satisfies the following
condition. \par

(A$\alpha$) There exists $a > 0$ and $A_{\alpha} > 0$ such that
\beq
\label{4.14}
v(r_1) - v(r_2) \geq \alpha^{-1} \ A_{\alpha}(r_2^{\alpha} - r_1^{\alpha}) \quad \hbox{for} \  0 <
r_1 < r_2 \leq a \quad . \eeq
\vskip 5 truemm

Note that as soon as $V \in L^p$ for some $p < \infty$, (H3) implies that $V$ is nonnegative and
tends to zero at infinity, so that (H3) implies (H2) with $V_- = 0$. We next discuss the last
condition in (H3). \par

The condition (A$\alpha$) can be formulated for any real $\alpha \not= 0$. In all cases it means
that $v$ is sufficiently decreasing in $(0, a]$. If $v \in {\cal C}^1((0,a])$, it is equivalent to
the fact that 
\beq
\label{4.15}
- v'(r) \geq A_{\alpha} \ r^{\alpha - 1} \qquad \hbox{for} \ 0 < r \leq a \quad .
\eeq
\noi It is easy to see that for any $\alpha$, $\beta$ with $\alpha \not= 0 \not= \beta$ and
$\alpha \leq \beta$, (A$\alpha$) implies (A$\beta$) with $A_{\beta} = a^{\alpha - \beta}
A_{\alpha}$. This is obvious from (\ref{4.15}) if $v \in {\cal C}^1((0, a])$. In the general case,
it reduces to the fact that $\alpha^{-1} (r_2^{\alpha} - r_1^{\alpha}) \geq \beta^{-1} a^{\alpha -
\beta} (r_2^{\beta} - r_1^{\beta})$ for $0 < r_1 \leq r_2 \leq a$ or equivalently, by making the
worst choice of $a$, namely $a = r_2$, and by scaling
$$\alpha^{-1}(1 - r^{\alpha}) \geq \beta^{-1}(1 - r^{\beta}) \qquad \hbox{for} \ 0 < r < 1 \quad ,
$$
\noi which can be verified easily. \par

>From the previous discussion it follows that the condition $\alpha \geq 2$ could be dropped in the
assumption (H3) without modifying that assumption. We have included that condition because we
shall use it in the subsequent applications. \par

For any real $\alpha \not= 0$, a potential $V(x) = C - \alpha^{-1} A_{\alpha} |x|^{\alpha}$ for
$0 < |x| \leq a$ satisfies, actually saturates the condition (A$\alpha$). The condition (A$\alpha$)
in (H3) means that $V$ is not too flat at the origin. For instance the potential $V(x) = 1 - \exp
(- 1/|x|)$ does not satisfy (H3), because there is no $\alpha$ for which it satisfies (A$\alpha$).
\par

In order to justify some computations to be made below, we shall need to approximate potentials
satisfying (H3) by potentials in ${\cal C}^1({I\hskip-1truemm R}^n\setminus \{0\})$ still
satisfying that condition. The standard approximation by convolution such as that used in the
proof of energy conservation is not suitable for that purpose because it does not preserve in an
obvious way the last condition of (H3). We proceed instead as follows. Let $\varphi \in {\cal
C}^{\infty}$ be radial nonnegative supported in the unit ball and satisfy $\int dx \ \varphi (x) =
1$. For $V \in L_{loc}^1({I\hskip-1truemm R}^n)$ and $j \geq 2$, we define
\beq
\label{4.16}
V_j(x) = j^n |x|^{-n} \int dy \ V(y) \ \varphi \left ( j|x|^{-1} (x - y) \right )
\eeq 
\noi or equivalently
\beq
\label{4.17}
V_j(x) = \int dz \ V(x - j^{-1} |x|z) \ \varphi (z) \quad .
\eeq
\noi Obviously $V_j \in {\cal C}^{\infty} ({I\hskip-1truemm R}^n \setminus \{ 0 \})$. The previous
regularisation satisfies the following properties. \\

\noi {\bf Lemma 4.3.} \par
{\it (1) Let $1 \leq p \leq \infty$ and $V \in L^p$. Then $V_j \in L^p$ and} 
\beq
\label{4.18}
\parallel V_j \parallel_p \ \leq (1 - j^{-1})^{-1/p} \ \parallel V \parallel_p \quad .
\eeq 
\noi {\it Let $1 \leq p < \infty$. Then $V_j$ tends to $V$ in $L^p$ norm when $j \to \infty$. \par
(2) If $V$ is radial, then $V_j$ is radial. If in addition $V$ is nonincreasing, then $V_j$ is
non\-in\-crea\-sing. If in addition $V$ satisfies the condition (A$\alpha$) for some $\alpha > 0$,
then $V_j$ also satisfies (A$\alpha$) with $a$ replaced by $a_j = a(1 + j^{-1})^{-1}$ and with
$A_{\alpha}$ replaced by $A_{\alpha j} = A_{\alpha} (1 - j^{-1})^{\alpha}$.} \\

\noi {\bf Proof.} {\bf Part (1)}. The case $p = \infty$ is obvious on (\ref{4.17}). For $p = 1$,
we estimate 
$$\parallel V_j \parallel_1 \ \leq \int dx \ dz \ |V(y)| \ \varphi (z) = \int dy \ dz |Dy/Dx|^{-1} \
|V(y)| \ \varphi (z)$$
\noi where $y = x - j^{-1}|x|z$ and $|Dy/Dx|$ is the Jacobian of the transformation $x \to y$ for
fixed $z$. One computes 
$$|Dy/Dx| = 1 - j^{-1} \widehat{x} \cdot z \geq 1 - j^{-1}$$
\noi which implies (\ref{4.18}) for $p = 1$. The general case of (\ref{4.18}) follows by
interpolation. \par

Convergence of $V_j$ to $V$ in $L^p$ for $p < \infty$ follows from the identity 
$$V_j(x) - V(x) = \int dz \ \varphi (z) \left ( V(x) - V(x - j^{-1} |x|y ) \right )$$
\noi for $V \in {\cal C}^1$ with compact support and follows from that special case for general $V$
by a density argument. \par

\noi {\bf Part (2).} Obviously $V$ radial implies $V_j$ radial. Let now $x_1$ and $x_2$ be
collinear, $x_1 = |x_1|\widehat{x}$, $x_2 = |x_2|\widehat{x}$ with $0 < |x_1| < |x_2|$. Then 
$$V_j(x_1) - V_j(x_2) = \int dz \ \varphi (z) \left (  V(x_1 - j^{-1} |x_1|z ) - V
( x_2 - j^{-1} |x_2|z ) \right )\quad .$$
\noi Clearly $x_i - j^{-1} |x_i|z = |x_i| (\widehat{x} - j^{-1}z)$, $i = 1,2$, are collinear and in
the same ratio as $x_1$ and $x_2$. Therefore $V_j$ is nonincreasing if $V$ is. Furthermore, if $V$
satisfies (\ref{4.14}) and if $|x_2|(1 + j^{-1}) \leq a$, then 
\begin{eqnarray*}
V_j(x_1) - V_j(x_2) &\geq & \alpha^{-1} A_{\alpha} \left ( |x_2|^{\alpha} - |x_1|^{\alpha} \right )
\int \varphi (z) \ dz \ |\widehat{x} - j^{-1} z |^{\alpha} \\
&\geq & \alpha^{-1} \ A_{\alpha} (1 - j^{-1})^{\alpha} \left ( |x_2|^{\alpha} - |x_1|^{\alpha} \right
) \end{eqnarray*}
\noi which completes the proof of Part (2). \par \nobreak
\hfill $\sq$ \\

In order to exploit the Morawetz inequality (\ref{4.6}), we shall need the following spaces. Let
$\sigma > 0$ and let $Q_i$ be the cube with edge $\sigma$ centred at $i\sigma$ where $i \in$
\Z$^n$ so that ${I\hskip-1truemm R}^n = \displaystyle{\mathrel{\mathop U_{i}}}
Q_i$. Let $1 \leq r$, $m \leq \infty$. We define 
$$\ell^m(L^r) = \left \{ u \in L_{loc}^r : \parallel u; \ell^m(L^r)\parallel \ = \ \parallel \ 
\parallel u; L^q(Q_i)\parallel ; \ell^m \parallel \ < \infty \right \} \quad .$$
\noi The space $\ell^m(L^r)$ does not depend on $\sigma$, and different values of $\sigma$ yield
equivalent norms. The previous spaces have been introduced by Birman and Solomjak \cite{4r}.
They allow for an independent characterization of local regularity and of decay at infinity in
terms of integrability properties. The H\"older and Young inequalities hold in those spaces, with
the exponents $m$ and $r$ treated independently. See \cite{11r} for more details. \par

We now turn to one of the main points of the proof, namely the fact that the Morawetz inequality
allows for the control of the $\ell^m(L^2)$ norm of the internal part of $u$ for some $m$. \\

\noi {\bf Proposition 4.2.} {\it Let $V$ satisfy (H1) with $n/4 < p_2 \leq p_1 < \infty$ and (H3)
and let $u \in X_{loc}^1({I\hskip-1truemm R})$ be a solution of the equation (\ref{1.1}). Then for
any $t_1$, $t_2 \in {I\hskip-1truemm R}$ with $1 \leq t_1 \leq t_2$, the following estimate holds 
\beq
\label{4.19}
\int_{t_1}^{t_2} dt (t\ {\rm Log} \ t + a )^{-1} \ \parallel u_<(t);\ell^{\alpha + 4}(L^2)
\parallel^{\alpha + 4} \ \leq C \ A_{\alpha}^{-1} \parallel u \parallel_2 \sqrt{E} \left ( \sqrt{E} +
\parallel u \parallel_2 \right )^{\alpha} \eeq

\noi where $u_<$ is defined by (\ref{4.2}) and where $C$ depends only on $n$, $\alpha$ and $a$.} \\

\noi {\bf Proof.} We first note that under the assumptions made on $V$, $u$ satisfies the
conservation laws of the $L^2$-norm and of the energy and belongs to $L^{\infty}({I\hskip-1truemm
R}, H^1)$ with $\parallel \nabla u ; L^{\infty}({I\hskip-1truemm R} , L^2) \parallel^2 \ \leq
2E$. The crucial steps in the proof can be easily performed if $V$ is ${\cal C}^1$
with compact support and we therefore approximate the actual $V$ by potentials of this type. We
first approximate $V$ by $V_j$ defined by (\ref{4.16}) and we then truncate $V_j$ both at the
origin and at infinity. Let $\psi_1 \in {\cal C}^{\infty}({I\hskip-1truemm R}^n)$, $\psi_1$ even,
with $0 \leq \psi_1 \leq 1$, $\psi_1 (x) = 1$ for $|x| \leq 1$ and $\psi_1 (x) = 0$ for $|x| \geq
2$ and define $\psi_{\ell} (x) = \psi_1 (x/\ell )$ for any $\ell > 0$. The final approximation
consists in approximating $V$ by $V_{jk} \equiv (1 - \psi_{1/k}) \psi_k \ V_j$ for large $j$ and
$k$. Clearly if $V \in L^p$ with $1 \leq p < \infty$, then $V_{jk}$ tends to $V_j$ in $L^p$ when $k
\to \infty$ while $V_j$ tends to $V$ in $L^p$ when $j \to \infty$ by Lemma 4.3 part (1). \par

We now define
\beq
\label{4.20}
J(V) = - \int_{t_1}^{t_2} dt \int dx \ \rho (x) \widehat{x} \cdot (V \star \nabla \rho )
\eeq
\noi so that by (\ref{4.6})
\beq
\label{4.21}
J(V) \leq 2 \parallel u \parallel_2 \sqrt{2E} \quad .
\eeq  
\noi We define $J(V_j)$ and $J(V_{jk})$ in the same way, for the same $\rho$. From estimates
similar to (\ref{4.10}) it follows that $J(V)$ is a continuous function of $V \in L^p$. In
particular $J(V_j) - J(V) \equiv \varepsilon_j$ tends to zero when $j \to \infty$. Clearly
\beq
\label{4.22}
J(V_j) \leq 2 \parallel u \parallel_2 \sqrt{2E} + \varepsilon_j \quad .
\eeq
\noi Similarly $J(V_{jk})$ tends to $J(V_j)$ when $k \to \infty$. Using (\ref{4.13}) and omitting
from now on the limits in the time integration, we rewrite $J(V_{jk})$ as
\bea
\label{4.23}
J(V_{jk}) = &-& \int dt \int dx \ \rho \ \widehat{x} \cdot \left ( (1 - \psi_{1/k}) \psi_k \ \nabla
V_j \star \rho \right ) \nn \\
&+& \int dt \int dx  \ \rho \ \widehat{x} \cdot \left ( (\nabla \psi_{1/k} - \nabla \psi_k )V_j 
\star \rho \right ) \quad .\eea
\noi We now estimate
\bea
\label{4.24}
\left | \int dx \ \rho \ \widehat{x} \cdot \left ( (\nabla \psi_{1/k}) \ V_j \star \rho \right )
\right | &\leq & C \parallel V_j \parallel_p \ \parallel u \parallel_s^4 \ \parallel \nabla
\psi_{1/k} \parallel_{\ell} \nn \\
&\leq & C \ k^{1-n/\ell} \ \parallel V_j\parallel_p \ \parallel u \parallel_s^4 \eea 
\noi with $4 \delta (s) = n/p + n/\ell$. We choose $5/4 < \delta (s) \leq 3/2$ which for $n/p \leq
4$ implies $n/\ell > 1$ and also ensures the local time integrability of the RHS of (\ref{4.24}),
so that the contribution of the LHS thereof to (\ref{4.23}) tends to zero when $k \to \infty$. The
contribution of $\nabla \psi_k$ to (\ref{4.23}) is treated in the same way. We obtain therefore 
\beq
\label{4.25}
J(V_j) = \lim_{k \to \infty} - \int dt \int dx \ \rho \ \widehat{x} \cdot \left ( (1- \psi_{1/k})
\psi_k \ \nabla V_j \star \rho \right ) \quad . \eeq
\noi We now use the fact that $V_j$ is radial with $V_j(x) \equiv v_j (|x|)$ and that $\psi_1$ is
even to rewrite (\ref{4.25}) as
\bea
\label{4.26}
J(V_j) = \lim_{k \to \infty} - {1 \over 2} \int dt \int dx \ dy \ \rho (x) \ \rho (y) (\widehat{x}
- \widehat{y}) \cdot (x- y) |x - y|^{-1}&& \nn \\
v'_j (|x - y|) \left ( (1 - \psi_{1/k} ) \psi_k \right ) (x - y)&&  . \eea 
\noi Now $v'_j \leq 0$ by (H3) and Lemma 4.3 part (2), while 
\beq
\label{4.27}
(\widehat{x} - \widehat{y}) \cdot (x - y) = \left ( |x|\  |y| - x \cdot y) (|x|^{-1} + |y|^{-1}
\right ) \geq 0 \eeq
\noi so that the integrand in (\ref{4.26}) is nonnegative and therefore nondecreasing in $k$.
Therefore, by the monotone convergence theorem 
\beq
\label{4.28}
J(V_j) = - (1/2) \int dt \int dx \ dy \ \rho (x) \ \rho (y) (\widehat{x} - \widehat{y}) \cdot (x -
y) |x - y|^{-1} v'_j(|x - y|) \quad . \eeq
\noi We now use the assumption (H3) applied to $V_j$ according to Lemma 4.3 part (2) and the
estimate (\ref{4.22}) to conclude that 
$$\int dt \int_{|x-y| \leq a_j} dx \ dy \ \rho (x) \ \rho (y) (\widehat{x} - \widehat{y})\cdot (x -
y) |x - y|^{\alpha - 2} \leq 2 A_{\alpha j}^{-1} (2\parallel u \parallel_2 \sqrt{2E} +
\varepsilon_j) \quad .$$
\noi Taking the limit $j \to \infty$ yields
\beq
\label{4.29}
\int dt \int_{|x - y| \leq a} dx \ dy \ \rho (x) \ \rho (y) (\widehat{x} - \widehat{y})\cdot (x -
y) |x - y|^{\alpha - 2} \leq 4 A_{\alpha}^{-1} \parallel u \parallel_2 \sqrt{2E} \quad . \eeq
\noi The next step in the proof consists in showing that the LHS of (\ref{4.29}) controls a
suitable $\ell^m(L^2)$ norm of $u_<$. Let $y_{\parallel}$ and $y_{\bot}$ be the components of $y$
parallel and normal to $x$. Then
\beq
\label{4.30}
{|x|\  |y| - x \cdot y \over |x|} = {|x| y_{\bot}^2 \over |x| \ |y| + x \cdot y} \geq {y_{\bot}^2
\over 2 |y|} \quad . \eeq
\noi Substituting (\ref{4.27}) and (\ref{4.30}) into (\ref{4.29}) and using the fact that $|y| \leq
|x| + a$ and $|y_{\bot}| \leq |x - y|$, we obtain 
\beq
\label{4.31}
\int dt \int dx \ \rho (x) (|x| + a )^{-1} \int_{|x - y| \leq a} dy \ \rho (y) |y_{\bot}|^{\alpha}
\leq 4 A_{\alpha}^{-1} \parallel u \parallel_2 \ \sqrt{2E} \quad . \eeq
\noi We now derive a lower bound of the integral over $y$ for fixed $x$. We first restrict that
integral to the cylinder $C_x$ of center $x$ and axis $x$ with diameter and height $a\sqrt{2}$.
For fixed $y_{\parallel}$, we consider the integral over $y_{\bot}$ which takes place in the ball
$B$ of radius $a/\sqrt{2} \equiv a_1$ centered at the origin in ${I\hskip-1truemm R}^{n-1}$. Let
$r_{\bot} = |y_{\bot}|$ and let $w$ be the vector field in ${I\hskip-1truemm R}^{n-1}$
$$w = y_{\bot} \left ( a_1^{\alpha} - r_{\bot}^{\alpha} \right )$$
\noi so that
$$\nabla \cdot w = (n - 1) a_1^{\alpha} - (n - 1 + \alpha )r_{\bot}^{\alpha} \quad .$$
\noi We now write
$$\int_B dy_{\bot} \ \rho \ \nabla \cdot w = - \int_B dy_{\bot} \ w \cdot \nabla \rho$$
\noi so that
\bea
&&(n - 1) a_1^{\alpha} \int_B dy_{\bot} \ \rho = (n - 1 + \alpha ) \int_B dy_{\bot} \
r_{\bot}^{\alpha} \ \rho \ - \ \int_B dy_{\bot} (a_1^{\alpha} - r_{\bot}^{\alpha} ) y_{\bot} \cdot
\nabla \rho \nn \\ &&\leq (n - 1 + \alpha ) \int_B dy_{\bot} \ r_{\bot}^{\alpha} \ \rho \ +
\ 2a_1^{\alpha} \parallel r_{\bot} u;L^2(B)\parallel \ \parallel \nabla u ; L^2(B)\parallel \quad .
\label{4.32} \eea \noi We next integrate (\ref{4.32}) over $y_{\parallel}$ and estimate the second
term in the RHS by applying the Schwarz inequality and extending the integral of $\nabla u$ to the
whole of ${I\hskip-1truemm R}^n$. We obtain
$$(n - 1) a_1^{\alpha} \parallel u; L^2 (C_x)\parallel^2 \ \leq (n - 1 + \alpha ) \parallel
r_{\bot}^{\alpha /2} u; L^2(C_x) \parallel^2 \ + \ 2a_1^{\alpha} \parallel \nabla u \parallel_2 \
\parallel r_{\bot} u; L^2(C_x) \parallel$$ 
$$\leq (n - 1 + \alpha ) \parallel r_{\bot}^{\alpha / 2} u;
L^2(C_x) \parallel^2 \ + \ 2a_1^{\alpha} \parallel \nabla u \parallel_2 \ \parallel r_{\bot}^{\alpha
/ 2} u; L^2(C_x) \parallel^{2/ \alpha} \parallel u ; L^2(C_x) \parallel^{1 - 2/\alpha}$$
\noi by the H\"older inequality, and therefore
$$
(n - 1) \parallel u; L^2(C_x)\parallel^{1 + 2/\alpha} \ \leq \ \parallel r_{\bot}^{\alpha / 2} u;
L^2(C_x)\parallel^{2/\alpha} \left \{ 2 \parallel \nabla u \parallel_2 + (n - 1 + \alpha )
a_1^{-1} \parallel u \parallel_2 \right \} \quad .$$
\noi Finally
\beq
\label{4.33}
\int_{C_x} dy \ r_{\bot}^{\alpha} \ \rho (y) = \ \parallel r_{\bot}^{\alpha / 2} u;
L^2(C_x)\parallel^2 \ \geq M \parallel u;L^2(C_x)\parallel^{\alpha + 2} \eeq
\noi with
\beq
\label{4.34}
M = (n - 1)^{\alpha} \left \{ 2 \sqrt{2E} + (n - 1 + \alpha ) \sqrt{2} a^{-1} \parallel u
\parallel_2 \right \}^{-\alpha} \quad . \eeq
\noi Substituting (\ref{4.33}) into (\ref{4.31}) yields
\beq
M \int dt \int dx \ \rho (x) (|x| + a)^{-1} \parallel u; L^2(C_x)\parallel^{\alpha + 2} \ \leq
4A_{\alpha}^{-1} \parallel u \parallel_2 \sqrt{2E} \quad .
\label{4.35}
\eeq
\noi We obtain a lower bound of the LHS of (\ref{4.35}) by replacing $u$ by $u_<$ and $(|x| +
a)^{-1}$ by $(t\ {\rm Log} \ t + a)^{-1}$ according to (\ref{4.2}). Introducing in addition the
decomposition of ${I\hskip-1truemm R}^n$ in unit cubes appropriate to the definition of
$\ell^m(L^2)$ spaces, we obtain
\beq
\label{4.36}
M \int dt (t \ {\rm Log} \ t + a)^{-1} \sum_i \int_{Q_i} dx |u_<(x)|^2 \parallel
u_<;L^2(C_x)\parallel^{\alpha + 2} \ \leq 4A_{\alpha}^{-1} \parallel u \parallel_2 \sqrt{2E} \quad .
\eeq 
\noi We next choose $\sigma$ in such a way that $C_x \supset Q_i$ for all $x \in Q_i$, and for that
purpose we take $\sigma = a(2n)^{-1/2}$, and we obtain 
\beq
\label{4.37}
M \int dt (t \ {\rm Log} \ t + a)^{-1} \ \parallel u_<;\ell^{\alpha + 4}(L^2) \parallel^{\alpha +
4} \ \leq 4A_{\alpha}^{- 1} \parallel u \parallel_2 \sqrt{2E} \quad . \eeq
\noi The estimate (\ref{4.19}) now follows from (\ref{4.37}) and (\ref{4.34}). \par \nobreak
\hfill $\sq$ \\

\noi {\bf Remark 4.1.} If the potential $V$ is flat at the origin, it seems difficult to extract a
norm estimate of $u$ directly from the inequality (\ref{4.6}). In fact if $V$ is ${\cal C}^1$ with
compact support, radial and nonincreasing with Supp $v' \subset \{r:0 < a \leq r \leq b \}$ and if 
\beq
\label{4.38}
{\rm Supp} \ u \subset \mathrel{\mathop \cup_{i \in \hbox{\Z}^n}} \ B ((a + b)i,
a/2) \eeq
\noi where $B(x, r)$ is the ball of center $x$ and radius $r$, then the RHS of (\ref{4.13}) is
zero. Therefore in order to get some information from (\ref{4.6}), one would have to use again
properties of the evolution, for instance the fact that the support property (\ref{4.38}) cannot
be preserved in time. \\

The basic estimate (\ref{4.19}) is not convenient for direct application to the integral equation,
and we now derive a more readily usable consequence thereof (cf. Lemma 5 in \cite{20r} and Lemma 5.3
in \cite{11r}). \\

\noi {\bf Lemma 4.4.} {\it Let $V$ satisfy (H1) with $n/4 < p_2 \leq p_1 < \infty$ and (H3) and let
$u \in X_{loc}^1({I\hskip-1truemm R})$ be a solution of the equation (\ref{1.1}). Then, for any
$t_1 \geq 1$, for any $\varepsilon > 0$ and for any $\ell \geq a$, there exists $t_2 \geq t_1 + \ell$
such that}
\beq
\label{4.39}
\int_{t_2- \ell}^{t_2} dt \parallel u_<(t);\ell^{\alpha + 4}(L^2)\parallel^{\alpha + 4} \ \leq
\varepsilon \quad . \eeq
\noi {\it One can find such a $t_2$ satisfying}
\beq
\label{4.40}
t_2 \leq \exp \left \{ \left ( 1 + {\rm Log}(t_1 + \ell ) \right ) \exp (M \ell / \varepsilon ) - 1
\right \} \eeq
\noi  {\it where $M$ is the RHS of (\ref{4.19}), namely}
\beq
\label{4.41}
M = C \ A_{\alpha}^{-1} \parallel u \parallel_2 \ \sqrt{E} \left ( \sqrt{E} \ + \parallel u
\parallel_2 \right )^{\alpha} \quad . \eeq
\vskip 5 truemm

\noi {\bf Proof.} Let $N$ be a positive integer. From (\ref{4.19}) we obtain
$$M \geq \sum_{j=1}^N \left [ (t_1 + j\ell ) \ {\rm Log} (t_1 + j\ell ) + a \right ]^{-1} \ K_j$$
\noi where
$$K_j = \int_{t_1 +(j -1)\ell}^{t_1 + j\ell} dt \parallel u_<(t);\ell^{\alpha +
4}(L^2)\parallel^{\alpha + 4} \quad .$$
\noi If $K_j \geq \varepsilon$ for $1 \leq j \leq N$, then
\begin{eqnarray*}
M &\geq& \varepsilon \sum_{j=1}^N \left [ \left ( t_1 + j\ell \right ) \ {\rm Log}\ \left ( t_1 +
j\ell \right ) + a \right ]^{-1} \\
&\geq & \varepsilon \ \ell^{-1} \int_{t_1 + \ell}^{t_1 + (N + 1)\ell} dt \left ( t \ {\rm Log} \ t +
a \right )^{-1} \\
& \geq&  \varepsilon \ \ell^{-1} \ {\rm Log} \ \left \{\left ( 1 + {\rm Log} \left (
t_1 + (N + 1)\ell \right ) \right ) \left ( 1 + {\rm Log} \left ( t_1 + \ell \right )
\right )^{-1} \right \}\end{eqnarray*}
\noi which is an upper bound on $N$, namely
$$t_1 + (N + 1) \ell \leq \exp \left \{ \left ( 1 + {\rm Log}\left ( t_1 + \ell \right ) \right )
\exp (M \ell / \varepsilon ) - 1 \right \} \quad .$$
\noi For the first $N$ not satisfying that estimate, there is a $j$ with $K_j \leq \varepsilon$ and
one can take $t_2 = t_1 + j\ell$ for that $j$. That $t_2$ is easily seen to satisfy (\ref{4.40}).
\par \nobreak
\hfill $\sq$ \\

We now exploit the estimates of Lemmas 4.2 and 4.4 together with the integral equation for $u$ to
prove that for $2 < r < 2^{\star}$, the $L^r$ norm of $u$ is small in large intervals. That part of
the proof follows the version given in \cite{7r} \cite{8r} for the NLS equation. \par

The assumptions made so far on $V$ include only (H1) with $n/4 < p_2 \leq p_1 < \infty$ and the
repulsivity condition (H3), but do not include any decay property at infinity. In order to proceed
further, we shall now require such properties, in the form of upper bounds on $p_1$. In particular
we require $p_1 < n$ for the next result (and subsequently $p_1 < n/2$). \\

\noi {\bf Lemma 4.5.} {\it Let $V$ satisfy (H1) with $1 \vee \ n/4 < p_2 \leq p_1 < n$ and (H3), and
let $u \in X_{loc}^1({I\hskip-1truemm R})$ be a solution of the equation (\ref{1.1}). Let $2 < r <
2^{\star}$. Then for any $\varepsilon > 0$ and for any $\ell_1 > 0$, there exists $t_2 \geq \ell_1$
such that} 
\beq
\label{4.42}
\parallel u ; L^{\infty}\left ( \left [ t_2 - \ell_1, t_2 \right ], L^r \right ) \parallel \ \leq
\varepsilon \quad . \eeq
\vskip 5 truemm

\noi {\bf Proof.} Since $u \in L^{\infty}({I\hskip-1truemm R}, H^1)$, it is sufficient to derive
the result for one value of $r$ with $2 < r < 2^{\star}$. The result for general $r$ will then
follow by interpolation with boundedness in $H^1$. We take one such $r$, with $0 < \delta \equiv
\delta (r) < 1$. Various compatible conditions will be imposed on $r$ in the course of the proof.
For technical reasons, we also introduce an $r_1 > 2^{\star}$, namely with $\delta_1 \equiv
\delta (r_1) > 1$, which will also have to satisfy various compatible conditions. \par

For future reference, we note that for any $s_1$, $s_2$, $t \in {I\hskip-1truemm R}$ 
$$G(s_1, s_2,t) = - i \int_{s_1}^{s_2} dt' \ U(t - t') \ f(u(t')) = U(t - s_2) \ u(s_2) - U(t - s_1)
\ u(s_1)$$ \noi so that
\beq
\label{4.43}
\parallel G(s_1, s_2, t)\parallel_2 \ \leq 2 \ \parallel u \parallel_2 \quad .
\eeq  

Let now $\varepsilon$ and $\ell_1$ be given. We introduce $\ell_2 \geq 1$, $t_1 > 0$ and $t_2 \geq
t_1 + \ell$ where $\ell = \ell_1 + \ell_2$, to be chosen later~; $\ell_2$ and $t_1$ will have to be
sufficiently large, depending on $\varepsilon$ but not on $\ell_1$ for given $u$. We write the
integral equation for $u$ with $t \in [t_2 - \ell_1, t_2]$ and we split the integral in that
equation as follows
\bea
\label{4.44}
u(t) &=& U(t) \ u(0) - i \left \{ \int_0^{t-\ell_2} dt' + \int_{t - \ell_2}^{t-1} dt' +
\int_{t-1}^t dt' \right \} U(t - t') \ f(u(t')) \nn \\
&\equiv& u^{(0)} + u_1 + u_2 + u_3
\eea
\noi in obvious notation, and we estimate the various terms in $L^r$ successively. Those estimates
will require auxiliary estimates of $f(u)$ in various spaces, and the latter will be postponed to
the end of the proof. In all the proof, $M$ denotes various constants depending only on $r$, $r_1$
and $\parallel u \parallel_2$, $E(u)$, possibly varying from one estimate to the next. \par

\noi {\bf Estimate of ${\bf u}^{\bf (0)}$.} \par \nobreak
It follows from Lemma 4.1 that
\beq
\label{4.45}
\varepsilon_0(t) \ \equiv \ \parallel U(t) \ u(0);L^{\infty}([t, \infty ), L^r)\parallel \to 0
\eeq  
\noi when $t \to \infty$, so that for $t \geq t_2 - \ell_1$
\beq
\label{4.46}
\parallel u^{(0)}(t)\parallel_r \ \leq \varepsilon_0(t) \leq \varepsilon_0 \left ( t_2 - \ell_1
\right ) \leq \varepsilon_0(\ell_2) \leq \varepsilon /4  \eeq
\noi for $\ell_2$ sufficiently large depending on $\varepsilon$. \par

\noi {\bf Estimate of ${\bf u}_{\bf 1}$.} \par \nobreak
We estimate by the H\"older inequality
\beq
\label{4.47}
\parallel u_1 \parallel_r \ \leq \ \parallel u_1 \parallel_2^{1 - \delta / \delta_1} \ \parallel
u_1 \parallel_{r_1}^{\delta / \delta_1} \quad . \eeq
\noi The $L^2$ norm of $u_1$ is estimated by (\ref{4.43}). The $L^{r_1}$ norm is estimated by the
use of the pointwise estimate (\ref{2.4}) as
\bea
\label{4.48}
\parallel u_1 \parallel_{r_1} &\leq& C \int_0^{t- \ell_2} (t - t')^{-\delta_1} \ \parallel f(u(t'))
\parallel_{\bar{r}_1} \nn \\
&\leq & C \left ( \delta_1 - 1 \right )^{-1} \ \ell_2^{1 - \delta_1} \ \parallel f(u); L^{\infty}({I\hskip-1truemm
R}, L^{\bar{r}_1} ) \parallel \quad . \eea
\noi We now use the estimate
\beq
\label{4.49}
\parallel f(u); L^{\infty}({I\hskip-1truemm R} , L^{\bar{r}_1})\parallel \ \leq M
\eeq 
\noi the proof of which is postponed. By (\ref{4.47}) (\ref{4.48}) (\ref{4.49}) we can ensure that
for $t_2 - \ell_1 \leq t \leq t_2$
\beq
\label{4.50}
\parallel u_1(t)\parallel_r \ \leq M \ \ell_2^{\delta / \delta_1 - \delta} \leq \varepsilon /4
\eeq
\noi for $\ell_2$ sufficiently large depending on $\varepsilon$. We now choose $\ell_2 = \ell_2
(\varepsilon )$ so as to ensure both (\ref{4.46}) and (\ref{4.50}). \par

We now turn to the estimate of $u_2$ and $u_3$. Here we need to consider the contributions of the
internal and external regions separately. We define
\beq
\label{4.51}
f_{<\atop >} (u) = u \left ( V \star |u_{<\atop >}|^2 \right )
\eeq
\noi and for $i = 2, 3$

\beq
\label{4.52}
u_i^{<\atop >}(t) = - i \int dt' \ U(t - t') \ f_{<\atop >}(u(t'))
\eeq
\noi where the time integral is performed in the appropriate interval, so that $f(u) = f_>(u) +
f_<(u)$ and therefore $u_i = u_i^> + u_i^<$. Note that this decomposition is not that defined by
(\ref{4.2}). \par

\noi {\bf Estimate of ${\bf u}_{\bf 2}$ and ${\bf u}_{\bf 3}$.} \par \nobreak
We estimate again by the H\"older inequality and (\ref{4.43})
\bea
\label{4.53}
\parallel u_2 \parallel_r  &\leq&  \parallel u_2 \parallel_2^{1 - \delta /\delta_1} \ \parallel u_2
\parallel_{r_1}^{\delta / \delta_1} \nn \\
&\leq & 2 \parallel u \parallel_2^{1 - \delta / \delta_1} \left ( \parallel u_2^>
\parallel_{r_1}^{\delta / \delta_1} \ + \ \parallel u_2^< \parallel_{r_1}^{\delta / \delta_1}
\right )  \eea
\noi and more simply
\beq
\label{4.54}
\parallel u_3 \parallel_r \ \leq \ \parallel u_3^> \parallel_r \ + \ \parallel u_3^< \parallel_r
\quad . \eeq

\noi {\bf Contribution of the external region.} \par \nobreak
By the same computation as in (\ref{4.48}), we estimate 
\beq
\label{4.55}
\parallel u_2^>(t)\parallel_{r_1} \ \leq C \left ( \delta_1 - 1 \right )^{-1} \ \parallel
f_>(u);L^{\infty} \left ( \left [ t - \ell_2, t \right ], L^{\bar{r}_1} \right ) \parallel \eeq
\noi and similarly
\beq
\label{4.56}
\parallel u_3^>(t) \parallel_r \ \leq C (1 - \delta )^{-1} \ \parallel f_>(u);L^{\infty} \left (
\left [ t - \ell_2 , t \right ], L^{\bar{r}} \right ) \parallel \quad . \eeq
\noi We next use the estimate
\beq
\label{4.57}
\parallel f_>(u(t));L^{\bar{r}_1} \cap L^{\bar{r}} \parallel \ \leq M \parallel
u_>(t)\parallel_2^{\mu} \eeq
\noi valid for some $\mu > 0$ and all $t \geq 1$, the proof of which is postponed. Using
(\ref{4.55})-(\ref{4.57}) and the propagation result of Lemma 4.2, part (2), we can ensure that the
contribution of the external region to $\parallel u_2(t)\parallel_r + \parallel u_3(t)\parallel_r$
which can be read on (\ref{4.53}) (\ref{4.54}) satisfies 
\beq
\label{4.58}
2 \parallel u \parallel_2^{1 - \delta / \delta_1} \ \parallel u_2^>(t)\parallel_{r_1}^{\delta /
\delta_1} + \parallel u_3^>(t) \parallel_r \  
\leq M \parallel u_>; L^{\infty} \left ( \left [ t_2 - \ell_1 - \ell_2 , t_2 \right ], L^2 \right
) \parallel^{\mu \delta /\delta_1} \leq \varepsilon /4 \eeq
\noi for all $t \in [t_2 - \ell_1, t_2]$ by taking $t_1$ sufficiently large, depending on
$\varepsilon$, since we have imposed $t_2 \geq t_1 + \ell_1 + \ell_2$. We now choose $t_1 = t_1
(\varepsilon )$ such that (\ref{4.58}) holds. \par

\noi {\bf Contribution of the internal region.} \par \nobreak

We shall use the estimate
\beq
\label{4.59}
\parallel f_<(u(t));L^{\bar{r}_1} \cap L^{\bar{r}}\parallel \ \leq M \parallel
u_<(t);\ell^m(L^2)\parallel^{m/s} \eeq
\noi where $m = \alpha + 4$, valid for some $s$ with $0 < s^{-1} < 1 - \delta$ and for all $t \geq
1$, the proof of which is postponed. Using again the pointwise estimate (\ref{2.4}), we estimate
\bea
\label{4.60}
\parallel u_2^<(t)\parallel_{r_1}  &\leq& C \int_{t-\ell_2}^{t-1} dt' (t - t')^{- \delta_1} \
\parallel f_<(u(t'))\parallel_{\bar{r}_1} \nn \\
&\leq & M(\delta_1 \bar{s} - 1)^{-1} \left \{ \int_{t - \ell_2}^{t-1} dt' \ \parallel
u_<(t');\ell^m(L^2)\parallel^m \right \}^{1/s} \eea
\noi by the H\"older inequality in time and (\ref{4.59}), and similarly
\bea
\label{4.61}
\parallel u_3^<(t)\parallel_r &\leq & C \int_{t-1}^t dt' (t - t')^{-\delta } \ \parallel
f_<(u(t'))\parallel_{\bar{r}} \nn \\
&\leq & M(1 - \delta \bar{s})^{-1} \left \{ \int_{t-1}^t dt' \ \parallel
u_<(t');\ell^m(L^2)\parallel^m \right \}^{1/s} \quad . \eea 

We now use (\ref{4.60}) (\ref{4.61}) and we apply Lemma 4.4 to conclude that there exists $t_2 \geq
t_1 + \ell \equiv t_1 + \ell_1 + \ell_2$ (remember that $\ell_2$ and $t_1$ are already chosen,
depending on $\varepsilon$) such that the contribution of the internal region to $\parallel
u_2(t)\parallel_r + \parallel u_3(t)\parallel_r$, which can be read on (\ref{4.53}) (\ref{4.54}),
satisfies
\beq
\label{4.62}
2 \parallel u \parallel_2^{1 - \delta / \delta_1} \ \parallel u_2^<(t)\parallel_{r_1}^{\delta /
\delta_1} + \parallel u_3^<(t) \parallel_r \ \leq  M \left \{ \int_{t_2 - \ell }^{t_2} dt' \
\parallel u_<(t');\ell^m(L^2)\parallel^m \right \}^{\delta /s\delta_1} \leq \varepsilon /4 \eeq
\noi for all $t \in [t_2 - \ell_1, t_2]$. \par

Collecting (\ref{4.46}) (\ref{4.50}) (\ref{4.58}) and (\ref{4.62}) and comparing with (\ref{4.44})
(\ref{4.53}) and (\ref{4.54}) yields (\ref{4.42}). \par

It remains to prove the estimates (\ref{4.49}) (\ref{4.57}) and (\ref{4.59}) on $f$. Estimates of a
quantity involving $\bar{r}_{(1)}$ mean that we want the estimates both for $\bar{r}$ and for
$\bar{r}_1$. \par

\noi {\bf Proof of (4.49) and (4.57).} \par \nobreak
We consider only (\ref{4.57}), of which (\ref{4.49}) is the special case obtained by replacing
$u_>$ by $u$ and taking $\mu = 0$. We estimate
\bea
\label{4.63}
\parallel f_>(u)\parallel_{\bar{r}_{(1)}} &\leq & C \parallel V \parallel_p \ \parallel u
\parallel_{r_2} \ \parallel u_> \parallel_{r_3}^2 \nn \\
&\leq & C \parallel V \parallel_p \ \parallel u \parallel_{r_2} \ \parallel u_> \parallel_2^{\mu} \
\parallel u \parallel_{2^{\star}}^{2 - \mu} \eea
\noi with $\delta_i = \delta (r_i)$, $i = 2, 3$, $0 \leq \delta_2 , \delta_3 \leq 1$, 
$$\delta_{(1)} + \delta_2 + 2 \delta_3 = n/p$$
\noi and $\mu = 2(1 - \delta_3)$. One can find admissible $r_2$ and $r_3$ provided $\delta_1 \leq
n/p \leq 3 + \delta - \mu$ which allows for $\delta < 1 < \delta_1$ and $\mu > 0$ provided $1 < n/p
< 4$. \par

\noi {\bf Proof of (4.59).} \par \nobreak
For radial nonincreasing $V$ satisfying (H1), one can decompose $V$ as $V = V_1 + V_2$ where $V_1
\in \ell^{p_1}(L^{\infty})$ and $V_2 \in \ell^1 (L^{p_2})$. One can take for instance $V_1(x) =
V(x) \ \chi (|x| \geq a)$ and $V_2(x) = V(x) \ \chi (|x| \leq a)$. Correspondingly, we decompose
$f_<(u) = f_{1<}(u) + f_{2<}(u)$. Using the fact that the spaces $\ell^k(L^r)$ are monotonically
increasing in $k$ and decreasing in $r$, we estimate
\bea
\label{4.64}
&&\parallel f_{1<}(u)\parallel_{\bar{r}_{(1)}} \ \leq C \parallel
f_{1<}(u);\ell^{\bar{r}_1}(L^2)\parallel \nn \\
&&\leq C \parallel V_1 ; \ell^{p_1}(L^{\infty})\parallel \ \parallel u\parallel_2 \ \parallel
u_<;\ell^k(L^2)\parallel^2 \eea
\noi provided
$$n/p_1 \geq \delta_{(1)} + 2 \delta (k) \quad .$$
\noi We then estimate
\beq
\label{4.65}
\parallel u_<;\ell^k(L^2)\parallel \ \leq C \parallel u \parallel_2^{1 - \lambda} \ \parallel
u_<;\ell^m(L^2)\parallel^{\lambda} \eeq
\noi for any $\lambda$ with
$$0 \leq \lambda \leq 1 \wedge \delta (k)/ \delta (m)$$
\noi so that $f_{1<}$ satisfies (\ref{4.59}) provided one can take $\delta (k) > 0$, namely 
provided $n/p_1 > \delta_1$. \par

We next estimate $f_{2<}$ as 
\beq
\label{4.66}
\parallel f_{2<}(u)\parallel_{\bar{r}_1} \ \leq C \parallel V_2;\ell^1(L^{p_2})\parallel \ \parallel
u \parallel_{r_2} \ \parallel u_<;\ell^q(L^{r_3})\parallel^2 \eeq
\noi with $0 \leq \delta_2$, $\delta_3 \leq 1$,
\beq
\label{4.67}
\delta_{(1)} + \delta_2 + 2 \delta_3 = n/p_2
\eeq
\beq
\label{4.68}
\delta_{(1)} + \delta_2 + 2\delta (q) = n \quad . 
\eeq
\noi If $\delta_2$ and $\delta_3$ satisfy (\ref{4.67}), and if $p_2 > 1$ one can use (\ref{4.68})
to define $q$ satisfying $r_3 < q < \infty$. If in addition $\delta_3 < 1$, one then estimates 
\beq
\label{4.69}
\parallel u_<; \ell^q(L^{r_3})\parallel \ \leq C \parallel u; \ell^k(L^2)\parallel^{1 - \delta_3} \
\parallel u \parallel_{2^{\star}}^{\delta_3} \eeq
\noi with $\delta (k) = (1 - \delta_3)^{-1}(\delta (q) - \delta_3)$, which together with
(\ref{4.65}) (\ref{4.66}) proves (\ref{4.59}) for $f_{2<}$. \par

It remains only to ensure (\ref{4.67}) with $\delta_3 < 1$. For that purpose we choose
$\delta_{(1)} + \delta_2 = n/p_2$ and $\delta_3 = 0$ if $\delta_{(1)} \leq n/p_2 \leq \delta_{(1)}
+ 1$, and $\delta_{(1)} + 1 + 2 \delta_3 = n/p_2$ if $\delta_{(1)} + 1 \leq n/p_2$, which ensures
$\delta_3 < 1$ provided $n/p_2 < 3 + \delta$. \par

Finally the required estimates (\ref{4.49}) (\ref{4.57}) and (\ref{4.59}) hold provided
$$\delta_1 < n/p_1 \leq n/p_2 < n \wedge (3 + \delta )$$
\noi which can be ensured under the assumptions made on $V$ by taking $\delta$ and $\delta_1$
sufficiently close to 1. \par \nobreak
\hfill $\sq$ \\

The next step in the proof consists in showing that the $L^r$ norm of $u(t)$ tends to zero when $t
\to \infty$. For that purpose we need to reinforce the decay assumption on $V$ at infinity to its
final form, namely $p_1 < n/2$. \\

\noi {\bf Lemma 4.6.} {\it Let $V$ satisfy (H1) with $ 1 \vee n/4 < p_2 \leq p_1 < n/2$ and (H3)
and let $u \in X_{loc}^1({I\hskip-1truemm R})$ be a solution of the equation (\ref{1.1}). Let $2 <
r < 2^{\star}$. Then $\parallel u(t)\parallel_r$ tends to zero when $t \to \infty$.} \\

\noi {\bf Proof.} The main step of the proof consists in showing that if $u$ satisfies (\ref{4.42})
for some $\varepsilon > 0$ sufficiently small (depending on $\parallel u \parallel_2$ and $E$) and
for some $\ell_1$ sufficiently large (depending on $u$ and on $\varepsilon$), then there exists $b$,
$0 < b \leq 1$, depending on $\varepsilon$ but not on $\ell_1$, such that 
\beq
\parallel u;L^{\infty}([t_2 - \ell_1, t_2 + b], L^r)\parallel \ \leq \varepsilon \quad .
\label{4.70}
\eeq

\noi For that purpose, we write the integral equation for $u$ with $t \in [t_2, t_2 + 1]$ and we
split the integral in that equation as follows 
\bea
\label{4.71}
u(t) &=& U(t) \ u(0) - i \left \{ \int_0^{t_2 - \ell_1} dt' + \int_{t_2 - \ell_1}^{t-1} dt' +
\int_{t-1}^{t_2} dt' + \int_{t_2}^t dt' \right \} U(t-t') \ f(u(t')) \nn \\
&\equiv & u^{(0)} + u_1 + u_2 + u_3 + u_4 \eea
\noi and we estimate the various terms successively. As in the proof of Lemma 4.5, auxiliary
estimates on $f(u)$ are postponed to the end. \par

\noi {\bf Estimate of ${\bf u}^{\bf (0)}$.} \par \nobreak
In the same way as in (\ref{4.46}), we can ensure that for $t \geq t_2$ 
\beq
\label{4.72}
\parallel u^{(0)}(t)\parallel_r \ \leq \varepsilon_0(t) \leq \varepsilon_0(t_2) \leq
\varepsilon_0(\ell_1) \leq \varepsilon /4
\eeq
\noi for $\ell_1$ sufficiently large depending on $u$ and on $\varepsilon$. \par

\noi {\bf Estimate of ${\bf u}_{\bf 1}$.} \par \nobreak
By the same estimates as in the proof of Lemma 4.5 (see especially (\ref{4.47}) (\ref{4.48}) with
$\ell_2$ replaced by $\ell_1 + t - t_2$, we can ensure that for $t \geq t_2$
\beq
\label{4.73}
\parallel u_1(t) \parallel_r \ \leq M \ \ell_1^{\delta / \delta_1 - \delta} \leq \varepsilon / 4
\eeq
 \noi for $\ell_1$ sufficiently large depending on $\varepsilon$, under the condition (\ref{4.49}).
\par

We require $\ell_1$ to be sufficiently large, depending on $\varepsilon$, to ensure (\ref{4.72}) and
(\ref{4.73}). \par

\noi {\bf Estimate of ${\bf u}_{\bf 4}$.} \par \nobreak
By the pointwise estimate (\ref{2.4}), we estimate
\bea
\label{4.74}
\parallel u_4(t) \parallel_r \ &\leq& C \int_{t_2}^t dt' (t - t')^{-\delta} \ \parallel f(u(t'))
\parallel_{\bar{r}} \nn \\
&\leq & C(1 - \delta )^{-1} (t- t_2)^{1 - \delta} \ \parallel f(u) ; L^{\infty} ({I\hskip-1truemm
R}, L^{\bar{r}}) \parallel \eea
\noi which by an estimate similar to (\ref{4.49}) with $\bar{r}_1$ replaced by $\bar{r}$, the proof
of which is omitted, enables us to ensure that 
\beq
\label{4.75}
\parallel u_4(t)\parallel_r \ \leq M \ b^{1 - \delta} \leq \varepsilon /4
\eeq
\noi for $t_2 \leq t \leq t_2 + b$ by taking $b$ sufficiently small depending on $\varepsilon$. \par

Note that the estimates made so far require only $p_1 < n$, but not the stronger condition $p_1 <
n/2$. The latter will be essential to estimate $u_2$ and $u_3$. \par

\noi {\bf Estimate of ${\bf u}_{\bf 2}$ and ${\bf u}_{\bf 3}$.} \par \nobreak
By the same method as in the proof of Lemma 4.5, especially (\ref{4.53}) (\ref{4.55}) (\ref{4.56})
we estimate 
\beq
\label{4.76}
\parallel u_2(t) \parallel_r \ \leq 2 \parallel u \parallel_2^{1 - \delta / \delta_1} \ \parallel u_2
(t) \parallel_{r_1}^{\delta /\delta_1}
\eeq
\beq
\label{4.77}
\parallel u_2(t)\parallel_{r_1} \ \leq C \left ( \delta_1 - 1 \right )^{-1} \ \parallel f(u) ;
L^{\infty} \left ( \left [ t_2 - \ell_1 , t_2 \right ], L^{\bar{r}_1} \right ) \parallel \eeq
\beq
\label{4.78}
\parallel u_3(t)\parallel_r \ \leq C (1 - \delta )^{-1} \parallel f(u) ;L^{\infty} \left ( \left [
t_2 - \ell_1 , t_2 \right ], L^{\bar{r}} \right ) \parallel \quad . \eeq
\noi We now use the estimates
\beq
\label{4.79}
\parallel f(u) \parallel_{\bar{r}_1}^{\delta / \delta_1} \ \leq M \parallel u \parallel_r^{1 + \nu}
\eeq
\beq
\label{4.80}
\parallel f(u) \parallel_{\bar{r}} \ \leq M \parallel u \parallel_r^{1 + \nu}
\eeq
\noi valid for some $\nu > 0$, the proof of which is postponed. \par

Using (\ref{4.76})-(\ref{4.80}) and the assumption (\ref{4.42}) made on $u$, we can ensure that
\bea
\label{4.81}
\parallel u_2 (t)\parallel_r \ + \ \parallel u_3(t)\parallel_r &\leq & M \parallel u;L^{\infty}
\left ( \left [ t_2 - \ell_1 , t_2 \right ], L^r \right ) \parallel^{1 + \nu} \nn \\
&\leq & M \ \varepsilon^{1 + \nu} \leq \varepsilon /4 \eea
\noi for all $t \in [t_2, t_2 + 1]$, for $\varepsilon$ sufficiently small depending only on
$\parallel u \parallel_2$ and $E$, namely for $\varepsilon \leq (4M)^{-1/\nu}$. \par

The proof of Lemma 4.6 now runs as follows. We pick an $\varepsilon$ satisfying the previous
condition. We next choose $\ell_1$ depending on $u$ and $\varepsilon$ so as to ensure (\ref{4.72})
(\ref{4.73}) for any $t_2 \geq \ell_1$, and we choose $b$ depending on $u$ and $\varepsilon$ so as
to ensure (\ref{4.75}). We next apply Lemma 4.5 with the previous $\varepsilon$ and $\ell_1$,
thereby obtaining $t_2$ such that (\ref{4.42}) holds. It then follows from (\ref{4.72}),
(\ref{4.73}), (\ref{4.75}) and (\ref{4.81}) that also (\ref{4.70}) holds. One then iterates the
argument with $\ell_1$ and $t_2$ replaced by $\ell_1 + jb$ and $t_2 + jb$ for $j = 1, 2, \cdots$,
which is possible since $\ell_1$ increases and $b$ is independent on $\ell_1$. One then obtains
$$\parallel u(t);L^{\infty} \left ( [ t_2 - \ell_1 , \infty ), L^r \right ) \parallel \ \leq
\varepsilon \quad .$$
\noi Applying that result for arbitrarily small $\varepsilon$ proves the Lemma. \par

It remains to prove the estimates (\ref{4.79}) (\ref{4.80}) for $f(u)$. For that purpose we
estimate
\beq
\label{4.82}
\parallel f(u) \parallel_{\bar{r}_{(1)}} \ \leq C \parallel V \parallel_p \ \parallel u
\parallel_{r_2}^3
\eeq
\noi where 
$$\delta_{(1)} + 3 \delta_2 = n/p \quad .$$
\noi Assuming without loss of generality that $n/p \leq 4 \delta$, we obtain $\delta_2 \leq
\delta$ in all cases, and we continue (\ref{4.82}) as
\beq
\label{4.83}
\cdots \leq C \parallel V \parallel_p \ \parallel u \parallel_2^{3(1 - \delta_2/\delta)} \
\parallel u \parallel_r^{3\delta_2/\delta} \eeq
\noi so that for both (\ref{4.79}) and (\ref{4.80}), the condition $\nu > 0$ becomes 
$$1 + \nu = 3 \delta_2 \ \delta_{(1)}^{-1} = \left ( n/p - \delta_{(1)} \right )
\delta_{(1)}^{-1} > 1$$
\noi or equivalently $n/p > 2 \delta_{(1)}$. Therefore (\ref{4.79}) (\ref{4.80}) hold provided $2
\delta_1 < n/p \leq 4 \delta$, which can be ensured under the assumptions made on $V$ by taking
$\delta$ and $\delta_1$ sufficiently close to 1. \par \nobreak
\hfill $\sq$ \\

We can now state the main result of this paper. \\

\noi {\bf Proposition 4.3.} {\it Let $V$ satisfy (H1) with $1 \vee n/4 < p_2 \leq p_1 < n/2$ and
(H3). \par
(1) Let $u \in X_{loc}^1({I\hskip-1truemm R})$ be a solution of the equation (\ref{1.1}). Then $u
\in X^1({I\hskip-1truemm R})$. \par
(2) The wave operators $\Omega_{\pm}$ are bijective bounded and continuous and their inverses
$\Omega_{\pm}^{-1}$ are bounded and continuous from $H^1$ to $H^1$.} \\

\noi {\bf Proof.} {\bf Part (1).} We give the proof in the special case where $V \in L^p$. The
general case of $V$ satisfying (H1) with $p_2 < p_1$ can be treated by a straightforward extension
of that proof, based on the fact that admissible pairs $(q_1, r_1)$ and $(q_2, r_2)$ in (\ref{2.8})
and (\ref{2.9}) are decoupled. Let $(q, r)$ be the admissible pair satisfying
$$1/2 < 2/q = \delta (r) = n/(4p) < 1 \quad .$$
\noi Let $0 < t_1 < t_2$. By the same estimates as in Sections 2 and 3 applied to the integral
equation for $u$ with initial time $t_1$, we estimate

\begin{eqnarray*}
y &\equiv & \parallel u; L^q([t_1, t_2], H_r^1) \parallel \ \leq C \parallel u(t_1);H^1 \parallel \\
&+& C \parallel V \parallel_p \ \parallel u;L^q([t_1, t_2], H_r^1) \parallel \ \parallel u; L^k
([t_1, t_2], L^r ) \parallel^2 \end{eqnarray*}
\noi with $1/k + 1/q = 1/2$, so that $k > q$ since $q < 4$. We interpolate 
$$\parallel u; L^k (L^r) \parallel^2 \ \leq \ \parallel u; L^q (L^r) \parallel^{2 - \lambda} \
\parallel u; L^{\infty} (L^r) \parallel^{\lambda}$$
\noi with $0 < \lambda = 4 - 8p/n < 2$, so that 
\beq
\label{4.84}
y \leq M + C \parallel V \parallel_p \ \parallel u; L^{\infty} ([t_1, t_2], L^r )
\parallel^{\lambda} \ y^{3 - \lambda} \quad . \eeq
\noi By Lemma 4.6, $\parallel u; L^{\infty}([t_1, t_2], L^r)\parallel$ can be made arbitrarily
small by taking $t_1$ sufficiently large, uniformly with respect to $t_2$. Furthermore for fixed
$t_1$, $y$ is a continuous (increasing) function of $t_2$, starting from zero for $t_2 = t_1$. It
then follows from (\ref{4.84}) that for $t_1$ sufficiently large $y$ is bounded uniformly in
$t_2$, namely that $u \in L^q([t_1, \infty ), H_r^1)$. Plugging that result again into the
integral equation yields that $u \in X^1({I\hskip-1truemm R}^+)$. The same argument holds for
negative times. \par

\noi {\bf Part (2).} The fact that $\Omega_+$ is a bijection of $H^1$ onto $H^1$ follows from the
fact that any initial data $u(0) = u_0 \in H^1$ generates a (unique) solution $u \in
X^1({I\hskip-1truemm R})$ by Part (1) of this proposition, so that $u$ has an asymptotic state $u_+ =
\displaystyle{\lim_{t \to \infty}} \widetilde{u}(t)$ by Proposition 3.2, part (1) and satisfies
the equation (\ref{3.3}) by Proposition 3.2 part (2) and therefore $u(0) = \Omega_+$ $u_+ \in {\cal
R}(\Omega_+)$. \par

Boundedness of $\Omega_+$ and $\Omega_+^{-1}$ follows from the conservation laws of the $L^2$
norm and of the energy, which together which (H1) (H3) imply that $\parallel u_0\parallel_2 \ = \ 
\parallel u_+\parallel_2$ and 
$$\parallel \nabla u_0 \parallel_2 \ \leq \ \parallel \nabla u_+ \parallel_2 \ = \sqrt{2E} \leq \
\parallel u_0;H^1 \parallel \ + C \parallel u_0 ; H^1 \parallel^2 \quad .$$

Continuity of $\Omega_+$ and $\Omega_+^{-1}$ follows from the corresponding properties in
Propositions 2.1 and 3.1. \par \nobreak
\hfill $\sq$ \\

\noi {\bf Remark 4.2.} It is an unfortunate feature of the method that it does not provide an
estimate of the norm of a solution $u$ in $X^1({I\hskip-1truemm R})$ in terms of the norm of
$u(0)$ in $H^1$, or equivalently of the norm of $u_+ = \widetilde{u}(\infty )$.

\end{document}